\documentclass[12pt]{article}
\usepackage{amsmath,amssymb,amscd}

\usepackage{amsmath}
\usepackage{graphicx}

\newtheorem{tm}{Theorem}[section]
\newtheorem{lm}[tm]{Lemma}
\newtheorem{co}[tm]{Corollary}
\newtheorem{re}[tm]{Remark}

\newtheorem{pr}[tm]{Proposition}
 
 \newenvironment{demo}[1]{\par\smallskip\par\begin{trivlist}
\item[]{\bf #1}\ }{\end{trivlist}\par\smallskip\par}
\newcommand{\Proof}{\begin{demo}{{\it Proof.\ }}}
\newcommand{\qed}{\end{demo}}
\newcommand{\toy}{\ \rule[0em]{0.5ex}{1.8ex}}
\newcommand{\QED}{\toy\end{demo}}
\newcommand{\la}{\langle}
\newcommand{\ra}{\rangle}
\newcommand{\nn}{\nonumber}
\newcommand{\III}{{\vert \kern-.10em \vert \kern-.10em \vert}}
\newcommand{\ve}{\varepsilon}
\newcommand{\al}{\alpha}

\makeatletter
 
 \@addtoreset{equation}{section}
\makeatother
 
\begin{document}
\setlength{\baselineskip}{15pt} 
%%%%%%%%%%%%%%%%%%%%%%%%%%%%%%%%%%%%%%%%%%%%%%%%%%%%%%%%%%%%%%%%%%%%%%
%
\bibliographystyle{plain}
\title{
Large deviation principle of Freidlin-Wentzell type 
for pinned diffusion processes
}
%%%%%%%%%%%%%%%%%%%%%%%%%%%%%%%%%%%%%%%%%%%%%%%%%%%%%%%%%%%%%%%%%%
%%%%%%%%%%%%%%%%%%%%%%%%%%%%%%%%%%%%%%%%%%%%%%%%%%%%%%%%%%%%%%%%%%%
\author{ 
Yuzuru INAHAMA
\footnote{
%e-mail:~\tt{inahama@math.nagoya-u.ac.jp} 
%\\
Graduate School of Mathematics,   Nagoya University
%\\
Furocho, Chikusa-ku, Nagoya 464-8602, JAPAN.
%\\
E-mail:~\tt{inahama@math.nagoya-u.ac.jp}
}
}
%\date{ \today }
\date{   }
% \pagestyle{empty}
%
% Start !!!
%
\maketitle
% \thispagestyle{empty}
%%%%\hspace{-5.5mm}

%%%%%%%%%%%%%%%%%%%%%%%%%%%%%%%%%%%%%%%%%%%%%%%%%%%%%%%%
%%%%%%%%%%%%%%%%%%%%%%%%%%%%%%%%
\begin{center}
{\bf  Abstract}
\end{center}
Since T. Lyons invented rough path theory, one of its most successful applications
is a new proof of Freidlin-Wentzell's large deviation principle for diffusion processes.
In this paper we extend this method to the case of pinned diffusion processes 
under a mild ellipticity assumption.
Besides rough path theory, our main tool is quasi-sure analysis, which is 
a kind of potential theory in Malliavin calculus.
%one of the deepest parts of Malliavin calculus.

%%%%%%%%%%%%%%%%%%%%%%%%%%%
%%%%%%%%%%%%%%%%%%%%%%%%%%%%%%%%%%%%%%%%%%%%%%%%%%%%%%%%%%%%%%%%%%
%%%%%%%%%%%%%%%%%%%%%%%%%%%%%%%%%%%%%%%%%%%%%%%%%%%%%%%%%%%%%%%
%%%%%%%%%%%%%%%%%%%%%%%%%%%%%%%%%%%%%%%%%%%%%%%%%%%%%%%%%%%%%%%%%%
%%%%%%%%%%%%%%%%%%%%%%%%%%%%%%%%%%%%%%%%%%%%%%%%%%%%%%%%%%%%%%%
%%%%%%%%%%%%%%%%%%%%%%%%%%%%%%%%%%%%%%%%%%%%%%%%%%%%%%%%%%%%%%%%%%
\section{ Introduction}
For the canonical realization of  $d$-dimensional Brownian motion $(w_t)_{0 \le t \le 1}$
and the vector fields 
$V_i : {\bf R}^n \to {\bf R}^n ~(1 \le i \le d)$ with sufficient regularity, 
let us consider the following  Stratonovich-type stochastic differential equation (SDE):
\[
dy_t =  \sum_{i=1}^d V_i (y_t)  \circ  d w_t^i 
\qquad
\mbox{with} \quad  y_0 =a \in {\bf R}^n.
\]
For the sake of simplicity, no drift term is added here, but modification is easy.
The correspondence $w \mapsto y$ is called the It\^o map 
and denoted by $y = \Phi (w)$.
It is well-known that the It\^o map is not continuous 
as a map from the Wiener space.
Moreover, it is not continuous with respect to
 any Banach norm on the Wiener space which preserves the structure of the Wiener space
 (see Sugita \cite{su2}).

Now, introduce a small positive parameter $\ve \in (0,1]$ and consider 
\[
dy_t^{\ve}  =  \sum_{i=1}^d V_i (y_t^{\ve})  \circ \ve d w_t^i 
\qquad
\mbox{with} \quad  y_0^{\ve} =a \in {\bf R}^n.
\]
%
%For simplicity of explanation, we do not add a drift term, but modification is easy.
%
Formally, $y^{\ve}  = \Phi (\ve w)$.
The process $(y_t^{\ve})_{0 \le t \le 1}$ takes its values in ${\bf R}^n$
and its law is a diffusion measure associated with the starting point $a$ and 
the generator ${\cal L}^{\ve} = (\ve^2 /2)\sum_{i=1}^d V_i^2$.

A classic result of Freidlin and Wentzell states 
the laws of $(y_t^{\ve})_{0 \le t \le 1}$ satisfies a large deviation principle 
as $\ve \searrow 0$.
The proof was not so easy. (See Friedman \cite{fried} 
or Dembo-Zeitouni \cite{dzbk} for instance).
If $\Phi$ were continuous, 
 we could use the contraction principle 
and the proof would be immediate from Schilder's large deviation principle for
the laws of $(\ve w_t)_{0 \le t \le 1}$.
However, it cannot be made continuous in the framework of the usual stochastic analysis.

Ten years ago, Ledoux, Qian, and Zhang \cite{lqz} gave a new proof by means of rough path theory,
which was invented by T. Lyons \cite{ly}. See also \cite{fvbk, lqbk} for example.
Roughly speaking, a rough path is a couple of a path itself and its iterated integrals.
Lyons established a theory of line integrals along rough paths
and ordinary differential equation (ODE) driven by rough paths.
The It\^o map in the rough path sense is deterministic and is sometimes called the Lyons-It\^o map.
The most important result in the rough path theory could be Lyons' continuity theorem
(also known as the universal limit theorem), which states that the Lyons-It\^o map is continuous 
in the rough path setting.
Brownian motion $(w_t)$ admits a natural lift to a random rough path $W$,
which is called Brownian rough path. 
If we put $W$ or $\ve W$ into the Lyons-It\^o map, 
then we obtain the solution of Storatonovich SDE $(y_t)$ or $(y^{\ve}_t)$, respectively.
%
%
%The most important result in the rough path theory could be  continuity theorem
%(also known as the universal limit theorem), which states that the Lyons-It\^o map is continuous 
%in the rough path setting.
%
%
They proved in \cite{lqz} that the laws of  $\ve W$ satisfy a large deviation principle of 
Schilder type with respect to the topology of the rough path space.
Large deviation principle of Freidlin-Wentzell type
for the laws of $(y^{\ve}_t)$ is immediate from this,
since the contraction principle can be used in this framework.
Since then many works on large deviation principle 
on rough path space have been published \cite{de, fv1, fv2, ik, mss}.

Then, a natural question arises;
can one obtain a similar result for pinned diffusion processes with this method, too?
More precisely, does the family of measures $\{ {\mathbb Q}^{\ve}_{a,a'}\}_{\ve >0}$
satisfy a large deviation principle as $\ve \searrow 0$?
Here, ${\mathbb Q}^{\ve}_{a,a'}$ is the pinned diffusion measure 
associated with ${\cal L}^{\ve}$, which starts at $a$ at time $t=0$
and ends at $a'$ at time $t=1$.
Heuristically, ${\mathbb Q}^{\ve}_{a,a'}$ is the law of $y_1^{\ve}$ 
under the conditional probability measure ${\mathbb P} (\,\cdot\, | y_1^{\ve}=a')$,
where ${\mathbb P}$ stands for the Wiener measure.
%
%(For a precise statement, see Theorem \ref{tm.ldp.main}.)

The aim of this paper is to answer this question affirmatively
under a certain mild ellipticity assumption for the coefficient vector fields.
Besides rough path theory, our main tool is quasi-sure analysis, which
is a sub-field of Malliavin calculus. 
It deals with objects such as Watanabe distributions (i.e., generalized Wiener functionals)
and capacities associated with Sobolev spaces.
Recall that motivation for developing this theory was 
to analyse (the pullbacks of) pinned diffusion measures on the Wiener space.
%
%
%It seems that our main result is not really new 
%and is an easy consequence of Gao and Ren's result \cite{gr}.
%This will be explained in the next section.

%%%%%%%%%%

%Takanobu and Watanabe presented this kind of large deviation principle
%under a hypoellipticity assumption for coefficient vector fields 
%(Theorem 2.1, \cite{tw}).
%
%This result seems very general and nice, but they gave no proof.
%Their tool are Malliavin calculus, and in particular, quasi-sure analysis.
%Recall that rough path theory did not exist, then.
%Presumably, they computed Besov norm of the solution of SDE, 
%but details are  unknown.

%Since we use rough path theory, we will compute, not the output, 
%but the input of the (Lyons-)It\^o map.
%Here, the input means $(w_t)$ itself and its iterated Storatonovich stochastic integrals.
%Hence, we believe that our proof via rough paths is probably simpler.
%Extending our method to the hypoellipticity case is an interesting 
%and important future task.

%%%%%%%%%%%

%Another preceding result is by  Hsu \cite{hsu} for a special case.
%He proved the case for (scaled) Brownian bridge on a complete Reimannian manifold $M$ 
%(i.e., the case 
%${\cal L}^{\ve} =(\ve^2 /2) \triangle_M$, where $\triangle_M$ stands for 
%Laplace-Beltrami operator on $M$).
%His proof is based on estimates and asymptotics for the heat kernel of $\triangle_M /2$ and
%no SDE appears in his paper.
%In this sense, this nice result of Hsu is not so "preceding" ours
%and it may not be very suitable to call it "Freidlin-Wentzell type".

%%%

The organization of this paper is as follows.
In Section 2, we give a precise setting, introduce assumptions,
 and state our main results. 
The most important among them is Theorem \ref{tm.ldp.sgt} ({\rm i}),
which is a large deviation principle on the geometric rough path space.
 We also discuss preceding results.
In Section 3,
we introduce Besov-type norms on the geometric rough path space
and prove their basic properties in 
relation to H\"older norms. 
In rough path theory, H\"older norms and variation norms are most important,
while Besov norms merely play auxiliary roles.
However, in connection with Malliavin calculus, Besov norms on the rough path space
become very important.
In Section 4,
we give a brief survey of Malliavin calculus and quasi-sure analysis for later use. 
We give basic facts on capacities and Watanabe distributions.
Sugita's theorem is the most important among them.
It states that a positive Watanabe distribution is actually a finite Borel measure
on the abstract Wiener space.

In Section 5, following Higuchi \cite{hi} and Aida \cite{ai},
we recall that Brownian motion $w$ admits a natural lift quasi-surely
and this version of Brownian rough path $W$ is $\infty$-quasi-continuous. 
In Section 6, 
using this quasi-continuity and Sugita's theorem, 
we obtain a probability measure on the geometric rough path space
such that its pushforward by the Lyons-It\^o map induces the pinned diffusion measure in question.

Section 7 is the the main part of our efforts. 
In this section  we prove our main results, in particular, Theorem \ref{tm.ldp.sgt} ({\rm i}).
%
%
%we consider the family  of finite Borel measures $\mu^{\ve}_{a,a'}$ on the geometric rough path space,
%which are the lift of the measures corresponding to positive Watanabe distributions
%$\delta_{a'} (y_1^{\ve}) = \delta_{a'} (y^{\ve} (1,a))$.
%
%which induces ${\mathbb Q}^{\ve}_{a,a'}$.
%
%In Theorem \ref{tm.ldp.sgt}, which is a key result in this paper
%and may be more important than our main result,
%we state a large deviation principle for $\{\mu^{\ve}_{a,a'} \}_{\ve >0}$.
%Notice that, in this theorem, ellipticity only at the starting point $a$ is assumed.
%
%
%Since the Lyons-It\^o map is continuous, the pushforward measures also 
%satisfy a large deviation principle.
%(We may consider the image measures of $\mu^{\ve}_{a,a'}$
% by another Lyons-It\^o map. See Remark \ref{re.anoI} below.)
%This is a special case of Takanobu and Watanabe's result (Theorem 2.1, \cite{tw}).
%
%Our main theorem (Theorem \ref{tm.ldp.main}) is almost immediate from this,
%when the pinned diffusion measures are well-defined and normalization is possible.
%
%
%Sections 8 and 9 are devoted for the proof of Theorem \ref{tm.ldp.sgt}.
%These two sections are the main part of our efforts.
%
%
Several key facts used in these sections are as follows;
{\rm (1)}~ large deviation estimate for capacities on geometric rough path space,
{\rm (2)}~ integration by parts formula in the sense of Malliavin calculus 
for Watanabe distributions, 
{\rm (3)}~ uniform non-degeneracy of Malliavin covariance matrix for solutions of the shifted SDE.

%%%%%%%%%%%%%%%%%%%%%%%%%%%%%%%%%%%%%%%%%%%%%%%%%%%%%%%%%%%%%%
%%%%%%%%%%%%%%%%%%%%%%%%%%%%%%%%%%%%%%%%%%%%%%%%%%%%%%%%%%%%%%%%%%
%%%%%%%%%%%%%%%%%%%%%%%%%%%%%%%%%%%%%%%%%%%%%%%%%%%%%%%%%%%%%%%
%%%%%%%%%%%%%%%%%%%%%%%%%%%%%%%%%%%%%%%%%%%%%%%%%%%%%%%%%%%%%%%%%%
%%%%%%%%%%%%%%%%%%%%%%%%%%%%%%%%%%%%%%%%%%%%%%%%%%%%%%%%%%%%%%%
%%%%%%%%%%%%%%%%%%%%%%%%%%%%%%%%%%%%%%%%%%%%%%%%%%%%%%%%%%%%%%%%%%
\section{ Setting and main result}

In this section we give a precise setting and state our main result.
Let $(w_t)_{0 \le t \le 1}$ be the canonical realization of 
$d$-dimensional Brownian motion.
We consider the following ${\bf R}^n$-valued Stratonovich-type SDE;
\begin{equation}\label{sde_set.eq}
dy_t^{\ve}  =  \sum_{i=1}^d V_i (y_t^{\ve})  \circ \ve d w_t^i + V_0 (\ve, y_t^{\ve}) dt
\qquad
\mbox{with} \quad  y_0^{\ve} =a \in {\bf R}^n.
\end{equation}
Here, $\ve \in [0,1]$ is a small parameter
and $V_i \in C_b^{\infty} ( {\bf R}^n, {\bf R}^n)$ for $1 \le i \le d$
and $V_0 \in C_b^{\infty} ( [0,1] \times {\bf R}^n, {\bf R}^n)$.
(A function is said to be of class $C_b^{\infty}$ if it is a bounded, smooth function 
with bounded derivatives of all order.)
For each $\ve$, $(y_t^{\ve})$ is a diffusion process 
with its generator
\[
{\cal L}^{\ve} = \frac{\ve^2}{2}  \sum_{i=1}^d V_i^2  + V_0 (\ve, \,\cdot\,).
\]
We denote by  ${\mathbb Q}^{\ve}_{a,a'}$ 
the pinned diffusion measure 
associated with ${\cal L}^{\ve}$
for $\ve>0$,  starting point $a$, and  terminal point $a'$ if it exists.

Let ${\cal H}$ be Cameron-Martin space for $(w_t)$. 
For $h \in {\cal H}$,
we denote by $\phi^0=\phi^0(h)$ be a unique solution of the following ODE;
\begin{equation}\label{ode_set.eq}
d \phi^0_t  =  \sum_{i=1}^d V_i (\phi_t^{0}) d h_t^i + V_0 (0, \phi_t^{0}) dt
\qquad
\mbox{with} \quad  \phi_0^{0} =a.
\end{equation}
We set ${\cal K}^{a,a'} = \{ h \in {\cal H} ~|~ \phi^0(h)_{1} =a' \}$.

Let $a, a'$ be given.
To exclude pathological cases, we assume the following  
\\
\\
{\bf (H0)}~ ${\cal K}^{a,a'}$ is not empty.
\\
\\
For the vector fields,
we assume  ellipticity only at the starting point $a$:
\\
\\
{\bf (H1)}~ The set of vectors 
$\{ V_1(a), \ldots, V_d(a) \}$ linearly spans ${\bf R}^n$.
\\
\\
Under this assumption, the heat kernel $p^{\ve}_t (a,a')$ exists 
for all $a' \in {\bf R}^n$, $t>0$ and $\ve>0$.
Note that $p^{\ve}_t (a,a')$ is continuous (actually smooth) in $a'$.
We also assume the following; 
% that the pinned diffusion measures are well-defined for given $a, a'$:
\\
\\
{\bf (H2)}~ $p^{\ve}_1 (a,a') >0$ for any $0 < \ve \le 1$.
\\
\\
Without {\bf (H2)} the normalization is impossible and 
 ${\mathbb Q}^{\ve}_{a,a'}$ or similar probability measures on the path space
cannot be defined. 
If the vector fields $\{ V_1, \ldots, V_d \}$ satisfy
everywhere ellipticity, then  {\bf (H0)}--{\bf (H2)} are satisfied.
%
%
%It is a unique probability measure on 
%$C_{a,a'} ([0,1], {\bf R}^n) :=\{ x \in  C([0,1], {\bf R}^n) ~|~ x_0 =a, x_1=a' \}$
%characterized by the following equation;
%for any $l \in {\bf N}_+$, any partition $\{ 0 =t_0 <t_1 <\cdots < t_l <t_{l+1} =1\}$ of $[0,1]$,
% and any bounded, continuous function $f:{\bf R}^{n \times l} \to {\bf R}$,
 %
 %
%\[
%\int f (x_{t_1}, \ldots, x_{t_l})  {\mathbb Q}^{\ve}_{a,a'}(dx)
%=
%p_1^{\ve} (a,a')^{-1}
%\int_{ ({\bf R}^n)^l }   f (a_1, \ldots, a_l) 
%\prod_{j=1}^{l+1}     p_{ t_j -t_{j-1} }^{\ve} (a_{j-1}, a_j)
% \prod_{j=1}^l da_j.
%\]
%Here, $a_0 =a$ and $a_{l+1} =a'$ by convention.
%
%
%In fact, ${\mathbb Q}^{\ve}_{a,a'}$ sits on 
%\[
%C^{\al -H}_{a,a'} ([0,1], {\bf R}^n) 
%=
%\{ x \in  C  ([0,1], {\bf R}^n) ~|~ \mbox{$x$ is $\al$-H\"older continuous, $x_0=a, x_1=a'$ } \}
%\]
%for any $\al \in (1/3, 1/2)$.
%
%
%Heuristically, ${\mathbb Q}^{\ve}_{a,a'}$ is the law of $y_1^{\ve}$ 
%under the conditional probability measure ${\mathbb P} (\,\cdot\, | y_1^{\ve}=a')$,
%where ${\mathbb P}$ stands for the Wiener measure.
%(This argument can be made rigorous with quasi-sure analysis, however.)

Before we state our main results,  let us introduce some notations.
We will give detailed explanations for 
terminologies from rough path theory and Malliavin calculus
in later sections.
Let $G\Omega^B_{\al, 4m} ({\bf R}^d)$ be the geometric rough path space 
with $(\al, 4m)$-Besov topology.
In this section we assume the following condition on the Besov parameters $(\al, 4m)$;
\begin{equation}\label{eq.amam}
m \in {\bf N}_+, \quad  \frac13 <\al < \frac12, \quad \al - \frac{1}{4m} > \frac13, \quad
\mbox{ and } \quad 4m -8m\al   >2.
\end{equation}
Note that $S_2 ({\cal H}) \subset G\Omega^B_{\al, 4m} ({\bf R}^d)$,
where $S_2$ stands for the natural (rough path) lift.
Under {\bf (H1)}, the positive Watanabe distribution $\delta_{a'} (y^{\ve}_1) =\delta_{a'} (y^{\ve}(1,a))$
is well-defined.
By Sugita's theorem, it is actually a finite Borel measure $\theta^{\ve}_{a,a'}$ on the Wiener space.
Its lift, denoted by $\mu^{\ve}_{a,a'} = (\ve S_2)_*\theta^{\ve}_{a,a'}$, 
is a finite Borel measure on $G\Omega^B_{\al, 4m} ({\bf R}^d)$.
Under  {\bf (H2)} the normalization of these measures are denoted by 
$\hat\theta^{\ve}_{a,a'}$ and $\hat\mu^{\ve}_{a,a'}$, respectively.

Set a rate function $I : G\Omega^B_{\al, 4m} ({\bf R}^n) \to [0, \infty]$ as follows;
\begin{align}
I(X) 
= 
\begin{cases}
    \|h\|^2_{{\cal H}}/2 & (\mbox{if $X= S_2(h)$ for some $h \in {\cal K}^{a,a'}$}), \\
    \infty &  (\mbox{otherwise}).
  \end{cases}
\nn
\end{align}
This rate function $I$ is good.
We also set $\hat{I}(X) = I(X) - \min\{ \|h\|^2_{{\cal H}}/2  ~|~h \in {\cal K}^{a,a'}\}$.
%Note here that the minimum is actually attained.

The following theorem is our main result in this paper.
The proof will be given in Section 7.
\begin{tm}\label{tm.ldp.sgt}
We keep the same notations and assumptions as above.
Let $m \in {\bf N}_+$ and  $1/3 <\al <1/2$ be as in (\ref{eq.amam}).
Then, the following {\rm (i)}--{\rm (ii)} hold:
\\
\noindent
{\rm (i)}~ 
 Assume {\bf (H0)} and {\bf (H1)}.
Then,
the family $\{ \mu^{\ve}_{a,a'}\}_{\ve >0}$ 
of finite measures
satisfies a large deviation principle 
on $G\Omega^B_{\al, 4m} ({\bf R}^d)$ as $\ve \searrow 0$ with a good rate function $I$, that is, 
for any Borel set $A \subset G\Omega^B_{\al, 4m} ({\bf R}^d)$,
the following inequalities hold;
\begin{align}
- \inf_{X \in A^{\circ} } I(X)  
\le
 \liminf_{\ve \searrow 0 } \ve^2 \log \mu^{\ve}_{a,a'} (A)
\le
 \limsup_{\ve \searrow 0 } \ve^2 \log \mu^{\ve}_{a,a'} (A)
\le 
- \inf_{X \in \bar{A} } I(X).
\nn
\end{align}
\noindent
{\rm (ii)}~ 
Assume {\bf (H0)}--{\bf (H2)}.
Then, the family $\{ \hat\mu^{\ve}_{a,a'}\}_{\ve >0}$ 
of probability measures 
satisfies a large deviation principle 
on $G\Omega^B_{\al, 4m} ({\bf R}^d)$ as $\ve \searrow 0$ with a good rate function $\hat{I}$.
\end{tm}

Note that Theorem \ref{tm.ldp.sgt} above also holds with respect to 
$\alpha'$-H\"older geometric rough path topology for any $\al' \in (1/3, 1/2)$,
since we can find $\al, m$ with (\ref{eq.amam}) 
such that $(\al, 4m)$-Besov topology is stronger than $\alpha'$-H\"older topology.

%\vspace{10mm}

Now we turn to the pinned diffusion measures.
Under {\bf (H1)} and {\bf (H2)} the probability measure $\hat\mu^{\ve}_{a,a'}$ exists.
Let $\Phi_{\ve}$ be the Lyons-It\^o map associated with $\{ V_1,\ldots, V_d; V_0(\ve, \,\cdot\,)  \}$.
$\tau_{\lambda}$ stands for the Young pairing map  with $\lambda_t =t$.
%Then, the probability measure $(a + \Phi_{\ve}^1)_* \hat\mu^{\ve}_{a,a'}$
%is well-defined.
%
%
Then, by the contraction principle and a slight generalization of Lyons' continuity theorem, 
$\{ (a + \Phi_{\ve}^1 \circ \tau_{\lambda})_* \hat\mu^{\ve}_{a,a'}\}_{\ve >0}$
satisfies a large deviation principle
on
\[
C^{\al -H}_{a,a'} ([0,1], {\bf R}^n) 
=
\{ x \in  C  ([0,1], {\bf R}^n) ~|~ \mbox{$x$ is $\al$-H\"older continuous, $x_0=a, x_1=a'$ } \}
\] 
for any $1/3 <\al <1/2$.
In fact, $(a + \Phi_{\ve}^1 \circ \tau_{\lambda})_* \hat\mu^{\ve}_{a,a'} 
= (\tilde{y}^{\ve})_* \hat\theta^{\ve}_{a,a'}$,
where $\tilde{y}^{\ve}$ is an $\infty$-quasi redefinition of $y^{\ve}$.

On the other hand, it is not easy to see  whether the pinned diffusion measure in an elementary sense 
exists under {\bf (H1)} and {\bf (H2)}.
But, for instance, 
if we assume everywhere hypoellipticity for $\{ V_1,\ldots, V_d; V_0(\ve, \,\cdot\,)  \}$
for any $\ve \in (0,1]$,
 then ${\mathbb Q}^{\ve}_{a,a'}$ exists and coincides with  
 $(a + \Phi_{\ve}^1 \circ \tau_{\lambda})_* \hat\mu^{\ve}_{a,a'}$.
In other words, $(a + \Phi_{\ve}^1 \circ \tau_{\lambda})_* \hat\mu^{\ve}_{a,a'}$
is a generalization of the pinned diffusion measure ${\mathbb Q}^{\ve}_{a,a'}$.
(See Section \ref{se.six} for details.)

Thus, as a corollary of Theorem \ref{tm.ldp.sgt}, 
we have obtained a Freidlin-Wentzell type large deviation principle 
for pinned diffusion measures.
Below we will state it explicitly. 
%
%(Note that assumptions of Theorem \ref{tm.ldp.main} implies {\bf (H2)}.)

%
Define  $J: C^{\al -H}_{a,a'} ([0,1], {\bf R}^n) \to [0, \infty]$ by 
$$
J(y) =\inf\{ \frac{\|h\|^2_{{\cal H}}  }{2} 
~|~ h \in {\cal K}^{a,a'}  \mbox{ with }  y = \phi^0 (h) \}
- \min\{ \frac{ \|h\|^2_{{\cal H}}  }{2} ~|~ h \in {\cal K}^{a,a'}  \}
$$ 
if $y =\phi^0 (h)$
for some $h \in {\cal K}^{a,a'}$
and 
define $J(y) = \infty$ if no such $h \in {\cal K}^{a,a'}$ exists.
$J$ is also good.
%Note that under {\bf (H1)} there is at most one $h \in {\cal K}^{a,a'}$ such that $y = \phi^0 (h)$
%for each given $y$.

%Now we state our main result in this paper.
%The rest of the paper will be devoted to proving this theorem.
%
%
\begin{tm}\label{tm.ldp.main}
Let $1/3 <\al <1/2$. In addition to {\bf (H0)}--{\bf (H2)}, assume 
everywhere hypoellipticity for $\{ V_1,\ldots, V_d; V_0(\ve, \,\cdot\,)  \}$ for any $\ve \in (0,1]$.
The family $\{ {\mathbb Q}^{\ve}_{a,a'}\}_{\ve >0}$ 
of pinned diffusion measures on $C^{\al -H}_{a,a'} ([0,1], {\bf R}^d)$
satisfies a large deviation principle 
as $\ve \searrow 0$ with a good rate function $J$, that is, 
for any Borel subset $A \subset C^{\al -H}_{a,a'} ([0,1], {\bf R}^n) $, 
\begin{align}
- \inf_{y \in A^{\circ} } J(y)  
\le
 \liminf_{\ve \searrow 0 } \ve^2 \log {\mathbb Q}^{\ve}_{a,a'} (A)
\le
 \limsup_{\ve \searrow 0 } \ve^2 \log {\mathbb Q}^{\ve}_{a,a'} (A)
\le 
- \inf_{y \in \bar{A} } J(y).
\nn
\end{align}
\end{tm}

%%%%%%%%%%%%%%%%%%%%%%%%%%%%
\begin{re}\label{re.anoI}
Let $F$ be any continuous map from $G\Omega^B_{\al,4m} ({\bf R}^{d})$ 
to a Hausdorff topological space.
Then, by the contraction principle,  
the family  $\{ F_* \mu^{\ve}_{a,a'} \}_{\ve >0}$ 
of image measures satisfies a large deviation principle as $\ve \searrow 0$.
Roughly speaking, 
if we take the Lyons-It\^o map determined by the vector fields $V_i$'s, 
then we get Theorem \ref{tm.ldp.main}.
But, we may take a different Lyons-It\^o map as $F$.  
This is worth mentioning since Theorem 2.1 \cite{tw} is formulated in such a way.
\end{re}

%%%%%%%%%%%%%%%%%%%%%%%%%%%%%%%%%%%%%%%%%%%%%%%%%%%%%%%%%%%%%%%

The rest of the section is devoted to looking at  preceding results.
It is clear that Theorem \ref{tm.ldp.sgt} is new,
because the number of paper on rough paths is not so large yet.
On the other hand, 
whether Theorem \ref{tm.ldp.main} is new is less clear
since pinned diffusions have been extensively studied.
However,  it seems to the author that
there is no explicitly written proof of this kind of large deviation.

%Using quasi-sure analysis,
%Gao and Ren \cite{gr} proved that the Freidlin-Wentzell estimate holds for 
%capacities on the Wiener space. 
%On the other hand, Sugita's theorem in \cite{su} states
%that any measure of finite energy (so, in particular, the conditional probability 
%given a non-degenerate Wiener functional) is controlled by a capacity. 
%Thus the Freidlin-Wentzell estimate holds for such measures. 
%
%
%By proving the "capacity version" of \cite{bs}, one can prove the equivalence of 
%Freidlin-Wentzell estimate and the large deviation principle.
%
%
% 
%Their result could remain true under the hypoellipticity assumption.

In 1994 
Takanobu and Watanabe presented this kind of large deviation principle
under a hypoellipticity assumption for coefficient vector fields 
(Theorem 2.1, \cite{tw}).
This result seems very general and nice, but they gave no proof.
Their tool are also quasi-sure analysis.
At that time rough path theory did not exist.
Presumably, they computed Besov norm of the solution of SDE, but details are unknown.
%
%
%(The author guesses that the results in \cite{gr} includes at least
%the most part of Theorem 2.1, \cite{tw}.
%But, 
%since there is no mention of \cite{tw} in \cite{gr}, 
%it is not so clear under what condidtion precisely these large deviations hold.)
%
(Extending our method to the hypoelliptic case is an interesting  and important future task.)

Another preceding result is by  Hsu \cite{hsu} for a special case.
He proved the case for (scaled) Brownian bridge on a complete Reimannian manifold $M$ 
(i.e., the case 
${\cal L}^{\ve} =(\ve^2 /2) \triangle_M$, where $\triangle_M$ stands for 
Laplace-Beltrami operator on $M$).
His proof is based on estimates and asymptotics for the heat kernel of $\triangle_M /2$ and
no SDE appears in his paper.
In this sense, this nice result of Hsu is not so "preceding" ours
and it may not be very suitable to call it "Freidlin-Wentzell type".
(After this work was done, the author was informed of Bailleul's result \cite{bai}.
Roughly speaking, this is a "hypoelliptic version" of Hsu's result \cite{hsu}.
His method is basically analytic, but he uses a little bit of rough path theory, too.)

A somewhat similar, but different result was obtained by Gao and Ren \cite{gr}.
They proved that the Freidlin-Wentzell large deviations hold for 
capacities on the Wiener space. 
It may be possible to prove that their result
implies large  deviations for pinned diffusions 
by showing the equivalence of the Freidlin-Wentzell large deviation 
and the Freidlin-Wentzell estimate. 
But, it is unclear how to carry it out
since there is no explanation in \cite{gr} or elsewhere.

Let us compare our method with the one in \cite{gr, tw}.
In \cite{gr, tw} they calculated with quasi-sure analysis 
the solution of SDE, which is the output of the It\^o map.
Since we use rough path theory, we will compute, not the output, 
but the input of the (Lyons-)It\^o map.
Here, the input means $(w_t)$ itself and its iterated Storatonovich stochastic integrals.
Hence, we believe that our proof via rough paths is probably simpler.
Extending our method to the hypoelliptic case is an interesting  and important future task.

On the other hand,  we do not use  the following;
{\rm (i)}~Malliavin calculus for Banach space-valued Wiener functionals as in \cite{gr, mn}
and 
{\rm (ii)}~Kusuoka-Stroock's large deviation transfer method for 
a family of measurable maps (Theorem 2.1, \cite{gr},  originally in Section 2, \cite{ks}).
Since both of them are hard to understand even for those who are familiar with
Malliavin calculus or large deviations,
being able to avoid them could be called an advantage.

%%%%%%%%%%%%%%%%%%%%%%%%%%%%%%%%%%%%%%%%%%%%%%%%%%%%%%%%%%%%%%%%%%
%%%%%%%%%%%%%%%%%%%%%%%%%%%%%%%%%%%%%%%%%%%%%%%%%%%%%%%%%%%%%%%
%%%%%%%%%%%%%%%%%%%%%%%%%%%%%%%%%%%%%%%%%%%%%%%%%%%%%%%%%%%%%%%%%%
%%%%%%%%%%%%%%%%%%%%%%%%%%%%%%%%%%%%%%%%%%%%%%%%%%%%%%%%%%%%%%%
%%%%%%%%%%%%%%%%%%%%%%%%%%%%%%%%%%%%%%%%%%%%%%%%%%%%%%%%%%%%%%%%%%
\section{Geometric rough path space with Besov norm}

In this section we introduce the geometric rough path space 
with Besov norm and recall its relation to the one with H\"older norm. 
Chapter 7 and Appendix A.2 in Friz and Victoir \cite{fvbk}
may be a nice summary on this issue.
In rough path theory, we usually use H\"older and variation norms.
In connection to Malliavin calculus, however, Besov norm play an essential role.

Throughout this paper, the time interval is $[0,1]$ and we set 
$\triangle =\{ (s,t) ~|~ 0 \le s \le t \le 1 \}$.
For any $Y \in C( \triangle,  {\bf R}^d)$,
we set 
\begin{align}
\|Y \|_{ \alpha -H}  &= \sup_{0 \le s < t \le 1  } \frac{| Y_{s,t} | }{ |t-s|^{\alpha}}
\qquad\qquad
&
(0< \alpha \le 1),
\label{hld.def}
\\
\| Y \|_{\alpha, m -B}  &= 
\Bigl(
\iint_{0 \le s < t \le 1 } 
 \frac{ | Y_{s,t} |^m }{  |t-s|^{1 +m\alpha}} dsdt  \Bigr)^{1/m}
\qquad\qquad
&
(m \ge 1, \, 0< \alpha \le 1).
\label{bsv.def}
\end{align}
These are called $\alpha$-H\"older norm and $(\alpha, m)$-Besov norm, respectively.
There are some variants of Besov-type norms, but we will use this one.

Set $C^{\alpha -H}([0,1], {\bf R}^d) =
 \{  x \in C([0,1], {\bf R}^d) ~|~ \|X^1 \|_{ \alpha -H} <\infty \}$,
where $X^1_{s,t} :=x_t -x_s$.
This is called the space of  ${\bf R}^d$-valued $\alpha$-H\"older continuous paths
and becomes a real Banach space with $\|x\| = |x_0| + \|X^1 \|_{ \alpha -H}$.
Its closed subset of paths that start at $a \in {\bf R}^d$ is denoted by 
$C_a^{\alpha -H}([0,1], {\bf R}^d)$.
In a similar way, if $\alpha - m^{-1} >0$, 
$C^{\alpha, m  -B}([0,1], {\bf R}^d) $ etc. are defined.
Then, $C^{\alpha, m -B}([0,1], {\bf R}^d)$
is continuously embedded in $C^{(\alpha-m^{-1}) -H}([0,1], {\bf R}^d)$. 
(See Appendix A.2, \cite{fvbk}.)
In this case,   closed subsets such as $C_a^{\alpha, m -B}([0,1], {\bf R}^d) $ are well-defined.

Let $T^2 ({\bf R}^d) = {\bf R} \oplus {\bf R}^d \oplus ({\bf R}^d)^{\otimes 2}$
be the truncated tensor algebra of step 2.
The set of elements of the form $(1, \bullet, \star)$ forms a non-abelian group 
under the tensor multiplication $\otimes$.
The unit element is ${\bf 1}=(1, 0,0)$.
Set 
$$
G^2 ({\bf R}^d)
=\{ (1, a_1, a_2) \in T^2 ({\bf R}^d) ~|~  a_2^{i,j}  +a_2^{j,i}   = a_1^i a_1^j \quad (1 \le i,j \le d)\}.
$$
It is easy to check that $G^2 ({\bf R}^d)$ becomes a subgroup.
It is called the free nilpotent group of step 2.
Note that the dilation on $T^2 ({\bf R}^d)$ 
(i.e., $(1, a_1, a_2) \mapsto (1, \lambda a_1, \lambda^2 a_2)$ for $\lambda \in {\bf R}$)
is well-defined on $G^2 ({\bf R}^d)$, too.

A continuous map
$X=(1, X^1, X^2) : \triangle \to   T^2 ({\bf R}^d)$
is called multiplicative if it satisfies that
\begin{align}
X^1_{s,t} = X^1_{s,u}+X^1_{u,t},  \qquad X^2_{s,t} = X^2_{s,u}+X^2_{u,t} +  X^1_{s,u} \otimes X^1_{u,t},
\qquad
\quad
(s \le u \le t).
\label{chen.eq}
\end{align}
This relation is called Chen's identity and can also be written as $X_{s,t} = X_{s,u} \otimes X_{u,t}$.
In particular, $X_{s,t}~(s \le t)$ is a "difference" of a group-valued path, since 
$X_{s,t} = (X_{0,s})^{-1} \otimes X_{0,t}$.

Let $1/3 <\alpha < 1/2$.
The space of ${\bf R}^d$-valued $\alpha$-H\"older rough path is defined by
\begin{align}
\Omega_{\alpha}^H ({\bf R}^d) 
&= \{  X=(1, X^1, X^2) \in C( \triangle, T^2( {\bf R}^d) )
\nn
\\ 
&  \qquad 
~|~  \mbox{multiplicative and   }  \|X^1 \|_{ \alpha -H} <\infty,  \|X^2 \|_{2 \alpha -H} <\infty\}.
\nn
\end{align}
The topology of this space is naturally induced by the following distance:
$d(X, Y) = \|X^1 - Y^1\|_{ \alpha -H} + \|X^2 -Y^2 \|_{2 \alpha -H}$.
In the same way, $(\alpha, m)$-Besov  rough path is 
defined for $m \ge 2$ and $1/3 <\alpha < 1/2$ with $\alpha -m^{-1}> 1/3$ as follows;
\begin{align}
\Omega_{m, \alpha}^B   ({\bf R}^d) 
&= \{  X=(1, X^1, X^2) \in C( \triangle, T^2( {\bf R}^d) )
\nn
\\ 
&  \qquad 
~|~  \mbox{multiplicative and   }  \|X^1 \|_{\alpha, m -B} <\infty, 
 \|X^2 \|_{ 2\alpha, m/2 -B} <\infty\}.
\nn
\end{align}
The topology of this space is naturally induced by the following distance:
$d(X, Y) = \|X^1 - Y^1\|_{\alpha, m -B} + \|X^2 -Y^2 \|_{  2 \alpha, m/2 -B}$.
In what follows, we will often write $X =(X^1, X^2)$ for simplicity, 
since the $0$th component "$1$" is obvious.

A Lipschitz continuous path (i.e., $1$-H\"older continuous path) 
$x \in C_0^{1 -H} ( [0,1], {\bf R}^d)$
admits a natural lift to a rough path by setting 
\[
X^1_{s,t} := x_t -x_s, \qquad X^2_{s,t} := \int_s^t (x_u -x_s) \otimes dx_u,  \qquad (s,t) \in \triangle. 
\]
It is easy to see that $X \in \Omega_{\alpha}^H ({\bf R}^d) \cap \Omega_{\alpha, m}^B ({\bf R}^d)$.
We call a rough path $X$ obtained in this way a smooth rough path lying above $x$,
or the lift of $x$.
The lift map is denoted by $S_2$, i.e., $X =S_2 (x)$.

An $\alpha$-H\"older weakly geometric rough path is $X \in \Omega_{\alpha}^H ({\bf R}^d)$
such that 
\begin{align}
X^{2; i,j}_{s,t} + X^{2; j, i}_{s,t} =  X^{1; i}_{s,t} X^{1; j}_{s,t},
\qquad
\qquad
(1 \le i,j \le d, \quad (s,t) \in \triangle).
\label{grpij.eq}
\end{align}
Here, $X^{2; i,j}_{s,t}$ stands for the $(i,j)$-component of $X^{2}_{s,t}$, etc.
Obviously, a smooth rough path satisfies (\ref{grpij.eq})
and is an $\alpha$-H\"older weakly geometric rough path.
The set of 
$\alpha$-H\"older weakly geometric rough paths is denoted by $G^w \Omega_{\alpha}^H ({\bf R}^d)$,
which is a closed subset of $\Omega_{\alpha}^H ({\bf R}^d)$.
In a similar way, $G^w \Omega_{\alpha, m}^B   ({\bf R}^d)$ is defined.
From (\ref{chen.eq}) and (\ref{grpij.eq}), a weakly geometric rough path is 
a "difference" of a $G^2 ({\bf R}^d)$-valued path.

Let $G\Omega_{\alpha}^H ({\bf R}^d)$ be the closure of the set of smooth rough paths, 
which is called the geometric rough path space with $\alpha$-H\"older norm.
The geometric rough path space $G \Omega_{ \alpha, m}^B   ({\bf R}^d) $ 
with $(\alpha, m)$-Besov norm is similarly defined.
Hence, we have the following inclusions;
\[
G\Omega_{\alpha}^H  ({\bf R}^d)  \subset G^w \Omega_{\alpha}^H   ({\bf R}^d) \subset 
\Omega_{ \alpha}^H   ({\bf R}^d),
\qquad
G\Omega_{\alpha, m}^B   ({\bf R}^d)  \subset G^w \Omega_{\alpha, m}^B   ({\bf R}^d) \subset 
\Omega_{\alpha, m}^B   ({\bf R}^d).
\]

%%%%%%
There is a natural correspondence between 
$G^2 ({\bf R}^d)$-valued path spaces in the usual sense and  weakly geometric rough path spaces. 
In order to explain that, we first introduce a distance $d$ on $G^2 ({\bf R}^d)$.
Let $\|{\bf a} \|$ be the Carnot-Caratheodory norm as in Chapter 7, \cite{fvbk}. 
This is a homogeneous norm on $G^2 ({\bf R}^d)$ with symmetry and subadditivity, too.
Its explicit form is not needed in this paper. 
For our purpose, it is enough to keep in mind that, there exist a constant $c>0$
such that 
\[
\frac{ 1}{c}\|{\bf a} \|  \le |a_1|_{{\bf R}^d} 
+ \sqrt{ |a_2|_{{\bf R}^{d} \otimes {\bf R}^{d}  }} \le c \|{\bf a} \|
\qquad
\mbox{for any ${\bf a} =  (1,a_1,a_2)  \in G^2 ({\bf R}^d)$}.
 \]
We set $d({\bf a}, {\bf b}) = \| {\bf a}^{-1} \otimes {\bf b}\|$.
This defines a left-invariant distance on $G^2 ({\bf R}^d)$, which induces the same topology 
as  the relative one inherited from $T^2 ({\bf R}^d)$.

%\begin{sloppypar}
%A $G^2 ({\bf R}^d)$-valued continuous path $g$ is said to be $\al$-H\"older continuous 
If $\sup_{0 \le s < t \le 1  } d(g_s, g_t) |t-s|^{-\alpha} <\infty$,
a $G^2 ({\bf R}^d)$-valued continuous path $g$ is said to be $\al$-H\"older continuous.
The set of such paths starting at ${\bf 1}$ is denoted by
$C^{\al -H}_{{\bf 1}} ([0,1], G^2 ({\bf R}^d))$.
Similarly, 
$g$ is said to be $(\al, m)$-Besov if 
$\iint_{0 \le s < t \le 1 } d(g_s, g_t)^m  |t-s|^{-( 1 +m\alpha)} dsdt <\infty$.
The set of such paths starting at ${\bf 1}$ is denoted by
$C^{\al, m -B}_{{\bf 1}} ([0,1], G^2 ({\bf R}^d))$.
For a weakly geometric rough path $X$, 
we can associate a $G^2 ({\bf R}^d)$-valued continuous path $t \mapsto X_{0,t}$.
This defines a natural bijection 
between $C^{\al -H}_{{\bf 1}} ([0,1], G^2 ({\bf R}^d))$ and $G^w \Omega_{\alpha}^H  ({\bf R}^d)$.
In the same way, there is a natural 
bijection 
between $C^{\al, m -B}_{{\bf 1}} ([0,1], G^2 ({\bf R}^d))$ and $G^w \Omega_{m, \alpha}^B  ({\bf R}^d)$.
Through these bijections, we introduce distance functions on 
$C^{\al -H}_{{\bf 1}} ([0,1], G^2 ({\bf R}^d))$ and 
$C^{m, \al -B}_{{\bf 1}} ([0,1], G^2 ({\bf R}^d))$, respectively.
%\end{sloppypar}

\begin{pr}\label{inj.HB.pr}
Assume $1/3< \al <1/2$, $m \ge 2$, and $\al -1/m >1/3$.
Then, $G\Omega_{\alpha, m}^B  ({\bf R}^d)$ is continuously embedded 
in $G\Omega_{\alpha -1/m}^H  ({\bf R}^d)$.
\end{pr}

\Proof
It is sufficient to show that 
$W\Omega_{\alpha, m}^B  ({\bf R}^d) \cong C^{\al, m -B}_{{\bf 1}} ([0,1], G^2 ({\bf R}^d))$
 is continuously embedded 
in $W\Omega_{\alpha -1/m}^H  ({\bf R}^d) \cong C^{(\al-m^{-1}) -H}_{{\bf 1}} ([0,1], G^2 ({\bf R}^d))$.
The inclusion,  which is called Besov-H\"older embedding, 
 is shown in Corollary A.2, \cite{fvbk}.
Its continuity is shown in Proposition A.9, \cite{fvbk}. 
\QED

It is easy to see from Proposition \ref{inj.HB.pr} above 
that $G\Omega_{\alpha, m}^B  ({\bf R}^d)$ is a complete separable metric space 
if $1/3< \al <1/2$, $m \ge 2$, and $\al -1/m >1/3$.

\begin{pr}\label{inj.cpt.pr}
Assume $1/3< \al < \al' <1/2$, $m \ge 2$, and $\al -1/m >1/3$.
Then, the injection 
$G\Omega_{\alpha', m}^B  ({\bf R}^d) \hookrightarrow G\Omega_{\alpha, m}^B  ({\bf R}^d)$
maps a bounded subset to a precompact subset. 
\end{pr}

\Proof
Let $\{ X(n)\}_{n =1}^{\infty}$ be any bounded sequence in $G\Omega_{\alpha', m}^B  ({\bf R}^d)$. 
Define $g(n)_t =  X(n)_{0,t}$.
Then, by Proposition \ref{inj.HB.pr},
$\{g(n) \}_{n =1}^{\infty}$ is a bounded sequence 
in $C^{(\al-m^{-1} )-H}_{{\bf 1}} ([0,1], G^2 ({\bf R}^d))$.
Hence, as functions of $t$, $\{g(n) \}_{n =1}^{\infty}$ is uniformly bounded 
and equicontinuous.
By Ascoli-Arzela's theorem, there exists a subsequence, which is denoted by $\{g(n) \}$ again, 
which converges to $g \in C_{{\bf 1}} ([0,1], G^2 ({\bf R}^d))$ in sup-norm.
Set $X_{s,t} = g_s^{-1} \otimes g_t$ for all $(s,t) \in \triangle$. 
It is sufficient to show that $\lim_{n \to \infty} X(n) =X$ in $(\al, m)$-Besov topology.
First, note that $\lim_{n \to \infty} \sup_{(s,t)  \in \triangle} |X(n)^j_{s,t} - X^j_{s,t}| =0$
for $j=1,2$.
Set $r =(1 +m\al')/(1+m\al) >1$ and let $r'$ be its conjugate exponent.
Then, by H\"older's inequality,
\begin{align}
\lefteqn{
\iint_{0 \le s< t \le 1 } 
 \frac{ | X(n)^1_{s,t}-  X(n')^1_{s,t} |^m }{  |t-s|^{1 +m\alpha}} dsdt  
 }
 \nn\\
 &\le
\sup_{(s,t)  \in \triangle} |X(n)^1_{s,t} - X(n')^1_{s,t}|^{m (1-1/r)}
\iint_{0 \le s< t \le 1 } 
 \frac{ | X(n)^1_{s,t}-  X(n')^1_{s,t} |^{m/r} }{  |t-s|^{1 +m\alpha}} dsdt 
 \nn\\
&\le
\sup_{(s,t)  \in \triangle} |X(n)^1_{s,t} - X(n')^1_{s,t}|^{m (1-1/r)}
\Bigl(
\iint_{0 \le s< t \le 1 } 
 \frac{ | X(n)^1_{s,t}-  X(n')^1_{s,t} |^{m} }{  |t-s|^{1 +m\alpha'}} dsdt 
 \Bigr)^{1/r}
 \frac{1}{2^{1/r'}}
  \nn\\
&\le
\mbox{const.} \times \sup_{(s,t)  \in \triangle} |X(n)^1_{s,t} - X(n')^1_{s,t}|^{m (1-1/r)}
\to 0
\qquad
\mbox{as $n, n' \to \infty$.}
\nn
 \end{align}
Here, we used that $(\al', m)$-Besov norm of $X(n)$ are bounded.
The limit in the Besov sense must coincide with $X^1$.
The second level path can be dealt with in the same way.
\QED

%%%%%%%%%%%%%%%%%%%%%%%%%%%%%%%%%%%%%%%%%%%%%%%%%%%%%%%%%%%%%%%%%%%%%%%%%%%%%%%%%%%%%%%%%%%%%%
%%%%%%%%%%%%%%%%%%%%%%%%%%%%%%%%%%%%%%%%%%%%%%%%%%%%%%%%%%%%%%%%%%%%%%%%%%%%%%%%%%%%%%
%%%%%%%%%%%%%%%%%%%%%%%%%%%%%%%%%%%%%%%%%%%%%%%%%%%%%%%%%%%%%%%%%%%%%%%%%%%%%%%%%%%%
%%  Section    quasi-sure Analysis
%%%%%%%%%%%%%%%%%%%%%%%%%%%%%%%%%%%%%%%%%%%%%%%%%%%%%%%%%%%%%%%%%%%%%%%%%%%%%%%%%%%%%%
%%%%%%%%%%%%%%%%%%%%%%%%%%%%%%%%%%%%%%%%%%%%%%%%%%%%%%%%%%%%%%%%%%%%%%%%%%%%%%%%%%%%%%
%%%%%%%%%%%%%%%%%%%%%%%%%%%%%%%%%%%%%%%%%%%%%%%%%%%%%%%%%%%%%%%%%%%%%%%%%%%%%%%%%%%%%%

\section{Preliminaries from quasi-sure analysis}

In this section we recall basics of Malliavin calculus and,
in particular, of quasi-sure analysis.
Generalized Wiener functionals (i.e., Watanabe distributions)
and capacities associated with Gaussian Sobolev spaces 
play important roles in this paper.
Analysis of these objects is called quasi-sure analysis.

\subsection{Basics of Malliavin calculus}
In this subsection,
we recall the basic notions in Malliavin calculus, in particular, the theory of Gaussian Sobolev spaces. 
We mainly follow Sections 5.8 --5.10, \cite{iwbk}, or Shigekawa \cite{sh}.

Let $({\cal W}, {\cal H}, \mu)$ be an abstract Wiener space.  
As usual, ${\cal H}$ and ${\cal H}^*$ are identified through Riesz isometry. 
Let $D$ be the ${\cal H}$-derivative and $D^*$ be its dual.
Ornstein-Uhlenbeck operator is denoted by $L =-D^* D$.
The first Wiener chaos associated with $h \in {\cal H}$ is denoted by $\la h,w \ra$.
If $\la h, \,\cdot\, \ra \in {\cal W}^*$, then it coincides with ${}_{{\cal W}^*}\la h,w \ra_{\cal W}$.
A function $F: {\cal W} \to {\bf R}$ is called a real-valued  polynomial if there exist 
$m \in {\bf N}$, $h_1, \ldots, h_m$, and  a real-valued polynomial $f$ of $m$-variables  
such that 
$F(w) = f( \la h_1 ,w \ra,   \ldots, \la h_m,w \ra )$.
The set of all real-valued polynomials are denoted by ${\cal P}$.
Let ${\cal K}$ be  a real separable Hilbert space.
$L^q ({\cal K}) =L^q ({\cal W}; {\cal K})$
denotes  the $L^q$-space of ${\cal K}$-valued functions.
  A function $G : {\cal W} \to {\cal K}$ is called a ${\cal K}$-valued  polynomial if there exist  
 $m \in {\bf N}_+$, $F_1, \ldots , F_m \in {\cal P}$, and $v_1, \ldots, v_m \in {\cal K}$ such that 
 $G(w)= \sum_{j=1}^m F_j (w) v_j$.
 The set of all ${\cal K}$-valued polynomials are denoted by ${\cal P} ({\cal K})$.

 For $q \in (1, \infty)$ and $r \in {\bf R}$, and $F \in {\cal P} ({\cal K})$, 
set $ \|  F\|_{q,r} =\|  (I -L)^{-r/2} F \|_q $.
 We define   the Sobolev space ${\bf D}_{q,r} ({\cal K})$ 
 to be the completion of ${\cal P} ({\cal K})$ with respect  to this norm.
When ${\cal K} ={\bf R}$, we simply write $L^q$, ${\bf D}_{q,r}$, etc.
If $q \le q'$ and $r \le r'$, then $\|\,\cdot \,\|_{q,r} \le \|  \,\cdot \, \|_{q',r'}$
and ${\bf D}_{q,r} ({\cal K}) \supset {\bf D}_{q',r'} ({\cal K})$.

%
%
%

%%%%%%%%%%%%%%%
When $r =k \in {\bf N}$, then the n Meyer equivalence holds (see \cite{sh} for example);
there exists a positive constant $c_{q,r}$ such that 
\[
c_{q,k}^{-1}  \|  F\|_{ {q,k}} 
\le 
\|F\|_q + \| D^k F\|_q
\le 
 c_{q,k}  \|  F \|_{ {q,k}}
\qquad
\mbox{for all $F\in {\bf D}_{q,k} ({\cal K})$.}
\]

%%%%%%%%%%%%%%   Wiener Chaos  %%%%%%%%%%%%%%%%

Now we discuss the Wiener chaos.
Set $\hat{\cal C}_n$ be the $L^2$-closure of 
real-valued polynomials of order $\le n$.
In particular, $\hat{\cal C}_0$ is the space of constant functions on ${\cal W}$.
$\hat{\cal C}_n$ is called the inhomogeneous Wiener chaos of order $n$.
Set ${\cal C}_0 =\hat{\cal C}_0$ and 
${\cal C}_n = \hat{\cal C}_n \cap {\cal C}_{n-1}^{\bot}$.
This is called the homogeneous Wiener chaos of order $n$.
Note that ${\cal C}_1 =\{ \la h, \,\cdot\, \ra ~|~ h \in {\cal H} \}$.
It is known that ${\cal C}_n$ is the eigenspace of $-L$ in $L^2$
that corresponds to the eigenvalue $n$. 
Clearly, $\hat{\cal C}_n = \oplus_{k=0}^n {\cal C}_k$.

\begin{pr}\label{chaos.nm.pr}
For any $n \in {\bf N}$, $q \in (1, \infty)$ and $r \ge 0$, there exists a constant $M_{n, q,r} \ge 1$
such that 
\begin{equation}\label{m.ch.ineq}
M_{n,q,r}^{-1} \|F\|_{2} \le  \|F\|_{q,r} \le  M_{n,q,r} \|F\|_{2}
\qquad
\qquad
(F \in \hat{\cal C}_n).
\end{equation}
In other words, restricted on $\hat{\cal C}_n$, all $(q,r)$-Sobolev norms are equivalent.
\end{pr}

\Proof
It is sufficient to show (\ref{m.ch.ineq}) for ${\cal C}_n$, instead of  $\hat{\cal C}_n$.
It is shown in Proposition 2.14, \cite{sh} that all $L^q$-norms are equivalent on ${\cal C}_n$.
Using the fact that ${\cal C}_n$ is an eigenspace of $-L$, we can easily see that 
all  $(q,r)$-Sobolev norms are equivalent on ${\cal C}_n$.
\QED

%%%%%%  Watanabe distributions %%%%%%%%%%%%

Now we recall  generalized Wiener functionals, which are also called Watanabe distributions.
Set ${\bf D}_{\infty} ({\cal K}) =\cap_{1<q<\infty, r \in {\bf R}} {\bf D}_{q,r} ({\cal K})$
and 
${\bf D}_{-\infty} ({\cal K}) =\cup_{1<q<\infty, r \in {\bf R}} {\bf D}_{q,r} ({\cal K})$.
Those are called the space of test functions and the space of generalized Wiener functionals, respectively.
We also use the following spaces of (generalized) Wiener functionals:
$\tilde{{\bf D}}_{\infty} ({\cal K}) =\cap_{k=1}^{\infty}  \cup_{1<q<\infty}   {\bf D}_{q,k} ({\cal K})$
and
$\tilde{{\bf D}}_{-\infty} ({\cal K}) =\cup_{k=1}^{\infty} \cap_{1<q<\infty}   {\bf D}_{q,-k} ({\cal K})$.
In other words,
$F \in \tilde{{\bf D}}_{\infty} ({\cal K})$ is equivalent to that, 
for any $k \ge 0$, there exists $q=q(k) >1$ such that $F \in {\bf D}_{q,k} ({\cal K})$.

Let $F =(F^1, \ldots, F^n) \in {\bf D}_{\infty} ({\bf R}^n) $.
The Malliavin covariance matrix is defined by $\bigl(  \la DF^i (w),DF^j (w)\ra_{{\cal H}} \bigr)_{1 \le i,j \le n}$.
We say that $F$ is non-degenerate in the sense of Malliavin if 
$\det \bigl(  \la DF^i (w),DF^j (w)\ra_{{\cal H}} \bigr)_{1 \le i,j \le n}^{-1}$ is in $L^q$ for all $1< q<\infty$.
This non-degeneracy is very important in Malliavin calculus.
For example,  such a non-degenerate $F$ can be composed with a Schwartz distribution $\phi$ defined on ${\bf R}^n$
and $\phi \circ F$ becomes a generalized Wiener functional.

We introduce seminorms on the space of tempered distributions on ${\bf R}^n$.
For $k \in {\bf Z}$ and 
a real-valued,
rapidly decreasing, smooth function $\phi$ of Schwartz class ${\cal S}({\bf R}^n)$ on ${\bf R}^n$, 
we define $\| \phi \|_{2k }=\| (1+|\,\cdot\,|^2 -\triangle /2 )^k \phi \|_{\infty}$.
Set ${\cal S}_{2k}({\bf R}^n)$ to be the completion of ${\cal S}({\bf R}^n)$ with respect to this norm.
Then ${\cal S}({\bf R}^n)= \cap_{k>0} {\cal S}_{2k}({\bf R}^n)$, which is a Fr\'echet space.
The dual space 
${\cal S}^{\prime}({\bf R}^n)= \cup_{k>0} {\cal S}_{-2k}({\bf R}^n)$ is called the Schwartz  space of tempered distributions.

S. Watanabe proved that, if a Wiener functional $F$ is non-degenerate, 
then the pullback map $\phi \mapsto \phi \circ F$ extends to a continuous linear  map between two distribution spaces.
The following is borrowed from pp. 378-379,  Ikeda and Watanabe \cite{iwbk}.

\begin{tm}\label{tm.pullback}
Let $F =(F^1, \ldots, F^n) \in {\bf D}_{\infty} ({\bf R}^n) $ be non-degenerate in the sense of Malliavin.
Then, for any $1< q <\infty$ and $k =0,1,2,\ldots$,
there exists a positive constant $C=C(q ,k,F)$ such that
\[
\|  \phi \circ F \|_{q, -2k}  \le C \| \phi \|_{-2k}
\qquad\qquad
(\phi \in {\cal S}({\bf R}^n))
\]
holds.
Therefore, the map $\phi \mapsto \phi \circ F$ extends  uniquely to a continuous linear map 
${\cal S}_{-2k}({\bf R}^n) \ni T  \mapsto T \circ F \in {\bf D}_{q,-2k}$.
In particular, $T \circ F \in  \cup_{k=1}^{\infty} \cap_{1< q <\infty} {\bf D}_{q,-2k} =\tilde{{\bf D}}_{- \infty}$.
\end{tm}

%%%%%%%%%%%%%

We recall an integration by parts formula in the context of Malliavin calculus 
(See p. 377, \cite{iwbk}).
For $F =(F^1, \ldots, F^n) \in {\bf D}_{\infty} ({\bf R}^n)$, we denote by 
$\tau^{ij} (w) =  \la DF^i (w),DF^j (w)\ra_{{\cal H}}$ the $(i,j)$-component of Malliavin covariance matrix.
We denote by $\gamma^{ij} (w)$ the $(i,j)$-component of the inverse matrix $\tau^{-1}$.
Note that $\tau^{ij} \in {\bf D}_{\infty} $ and
$D \gamma^{ij} = \sum_{k,l} \gamma^{ik} ( D\tau^{kl} ) \gamma^{lj} $.
Hence, derivatives of $\gamma^{ij}$ can be written in terms of
$\gamma^{ij}$'s and the derivatives of $\tau^{ij}$'s.
Suppose $G \in {\bf D}_{\infty}$ and $T \in {\cal S}^{\prime} ({\bf R}^n)$.
Then, the following integration by parts holds;
\begin{align}
{\mathbb E} \bigl[
(\partial_i T \circ F ) \cdot G 
\bigr]
=
{\mathbb E} \bigl[
(T \circ F ) \cdot \Phi_i (\, \cdot\, ;G)
\bigr]
\label{ipb1.eq}
\end{align}
where $\Phi_i (\, \cdot\, ;G) \in  {\bf D}_{\infty}$ is given by 
\begin{align}
\Phi_i (w ;G) &=
-\sum_{j=1}^d  
\Bigl\{
-\sum_{k,l =1}^d G(w) \gamma^{ik }(w)\gamma^{jl }(w)    \la D\tau^{kl} (w),DF^j (w)\ra_{{\cal H}}
\nn\\
&
\qquad\qquad
+
\gamma^{ij }(w) \la DG (w),DF^j (w)\ra_{{\cal H}} + \gamma^{ij }(w) G (w) LF^j (w)
\Bigr\}.
\label{ipb2.eq}
\end{align}
Note that the expectations in (\ref{ipb1.eq}) are in fact  
the pairing of $\tilde{{\bf D}}_{- \infty}$ and $\tilde{{\bf D}}_{\infty}$.

%%%%%%%%%%%%%%%%%%%%%%%%%%%%%%%%%%%%%%%%%%%%%
%%%%%%%%%%%%%%%%%%%%%%%%%%%%%%%%%%%%%%%%%%%%%

\subsection{Basics of capacities}

In this subsection, we recall the definition and basic properties of
the capacity associated with the Sobolev space ${\bf D}_{q,r}$
for $1<q <\infty$ and $r \in {\bf N}$. 

The contents of the subsection is borrowed from Chapter 9,  Malliavin \cite{ma}.
In this book, they work on a particular Gaussian space, namely ${\bf R}^{\infty}$ 
equipped with countable product of one-dimensional standard normal distribution.
But, the results in \cite{ma}, at least the ones we will use in this paper, 
hold true on any abstract Wiener space. 
In some literatures, a slightly different definition of capacities is used
(see \cite{su} for instance).

Let $({\cal W}, {\cal H}, \mu)$ be an abstract Wiener space.
We keep the same notation as in the previous subsection.
Throughout this subsection,  we set $1<q <\infty$ and $r \in {\bf N}$.  
For an open subset $O \subset {\cal W}$, we set 
\begin{equation}\label{capopen.def}
{\rm c}_{q,r} (O) 
=
\inf \{ \| \phi\|_{q,r}  ~|~ \mbox{ $\phi(w) \ge 1$ a.e. on $O$ and $\phi(w) \ge 0$ a.e. on ${\cal W}$} \}.
\end{equation}
For any subset $A \subset {\cal W}$, which is not necessarily open,
we define 
\begin{equation}\label{cap2open.def}
{\rm c}_{q,r} (A) 
=
\inf \{  {\rm c}_{q,r} (O) ~|~ \mbox{$A \subset O$ and  $O$ is open }  \}.
\end{equation}
Since $(q,r)$-norm is increasing in both $q$ and $r$, so is $c_{q,r} (A)$.
We say that a property $\pi= \pi_w$ holds $(q,r)$-quasi-everywhere 
if ${\rm c}_{q,r} (\{ w |\mbox{$\pi_w$ does not hold} \}) =0$.
We say that a property $\pi= \pi_w$ holds quasi-surely 
if it holds  $(q,r)$-quasi-everywhere for all $q \in (1, \infty)$ and $r \in {\bf N}$.
A subset $A$ is called slim if ${\rm c}_{q,r} (A) =0$ for all $q \in (1, \infty)$ and $r \in {\bf N}$.

Now let us discuss quasi-continuity.
A function $\phi$ from ${\cal W}$ to ${\bf R}$ (or to a  metric space)
is said to be $(q,r)$-quasi-continuous 
if, for any $\ve >0$, 
there exists an open set $O_{\ve}$ such that ${\rm c}_{q,r} (O_{\ve} ) <\ve$
and the restriction $\phi |_{ O^c_{\ve}}$ is continuous.
A function $\phi$ from ${\cal W}$ to ${\bf R}$ (or to a  metric space)
is said to be $\infty$-quasi-continuous 
if it is $(q,r)$-quasi-continuous for all $q$ and $r$.
This is equivalent to that
there is a decreasing sequence of open subsets $O_n$ such that 
$\lim_{n \to \infty} {\rm c}_{n,n} (O_n) =0$
and
$\phi |_{ O^c_{n}}$ is continuous.
(In this paragraph, the set $O_{\ve}$ or $O_n$ actually need not be open, 
since any subset can be approximated in terms of capacity
by an open subset from outside.)

For a measurable function $\psi$, 
we say $\psi^*$ is a $(q,r)$-redefinition  of $\psi$ if
$\psi^*$ is $(q,r)$-quasi-continuous and $\psi =\psi^*$ a.s.
($\psi^*$ is also called $(q,r)$-quasi continuous modification.)
Note that $(q,r)$-redefinition is essentially unique when it exists.
It is shown in Theorem 2.3.3, \cite{ma}  that 
any $\phi \in {\bf D}_{q,r}$ admits a $(q,r)$-redefinition.
Similarly, for a measurable function $\psi$, 
we say $\psi^*$ is a $\infty$-redefinition  of $\psi$ if
$\psi^*$ is $\infty$-quasi-continuous and $\psi =\psi^*$ a.s.
It is shown in Subsection 2.4, \cite{ma}  that 
any $\phi \in {\bf D}_{\infty}$ admits an (essentially unique) $\infty$-redefinition.

%%%%%%%%%%%%%%%

Now we give two useful basic lemmas without proofs for later use.
One is Borel-Cantelli's lemma (Corollary 1.2.4, \cite{ma}) and the other is Chebyshev's lemma (Theorem 2.2, \cite{ma}).
\begin{lm}\label{lm.borcan}
Let $1 <q <\infty$ and $r \in {\bf N}$.
Assume that $A_k \subset {\cal W}$  $(k \in {\bf N})$ satisfy that $\sum_k {\rm c}_{q,r} (A_k) <\infty$.
Then, 
${\rm c}_{q,r} (\limsup_{k \to \infty} A_k) =0$.
\end{lm}

\begin{lm}\label{lm.cheby}
For any $1 <q <\infty$ and $r \in {\bf N}$, there exists a positive constant $M_{q,r}$
such that, for any $\phi \in {\bf D}_{q,r}$ and any $R>0$,
we have
 \[
{\rm c}_{q,r} ( \{  w ~|~ \phi^* (w) >R\} )  \le R^{-1} M_{q,r} \| \phi \|_{q,r}.
\]
\end{lm}
%%%%%%%%%%%%%%%

In the finite dimensional calculus, it is known that
a positive Schwartz distribution is a measure.
An analogous fact is true in Malliavin calculus, too.
It is called Sugita's theorem (Theorem 3.0, \cite{ma} or Sugita \cite{su})  
and will play an important role in the sequel.

\begin{pr}\label{pr.sugita}
Let $l \in {\bf D}_{- \infty}$ be a positive Watanabe distribution,
that is, it satisfies 
${}_{ {\bf D}_{-\infty}} \la l, f \ra_{ {\bf D}_{\infty}}  ={\mathbb E}[l \cdot  f] \ge 0$ 
for any $f \in {\bf D}_{\infty}$ such that $f >0$ a.s.
Then, 
there exists a unique positive Borel measure $\theta$ of finite total mass  such that
\begin{equation}\label{sugita.eq}
\la l, g \ra  ={\mathbb E}[l \cdot  g] 
=
\int_{{\cal W}}  g^* (w) d\theta(w),
\qquad\qquad
(g \in {\bf D}_{\infty}).
\end{equation}
Furthermore, $\theta$ does not charge a slim set. 
(Hence, any choice of $\infty$-redefinition $g^*$ will do.)
\end{pr}
It is easy to see that,
if $l \in {\bf D}_{q', -r}$ with $1/q +1/q' =1$, 
then $\theta (O) \le {\rm c}_{q,r}(O) \| l \|_{q', -r}$ for any open set $O$.
(Note that $O$ need not be open here again.)

Let us discuss the equilibrium potential.
Given a Borel subset $A \subset {\cal W}$, we define
\begin{equation}\label{eq.libr.def}
{\cal F}_{q,r} (A)= 
\{  
u \in {\bf D}_{q,r} ~|~ \mbox{ $u^* \ge 1$ $(q,r)$-quasi-everywhere on $A$ }
\}.
\end{equation}
This is a closed convex subset of ${\bf D}_{q, r}$ and 
has a unique element $\phi_A \in {\cal F}_{q,r} (A)$, which minimizes $(q,r)$-Sobolv norm.
We call $\phi_A$ the equilibrium potential of $A$.
Then, Theorem 4.4, \cite{ma} states that, 
\begin{equation}\label{lib.pot.ineq}
{\rm c}_{q,r} (A) = \| \phi_A\|_{q,r}  
  \le \| u\|_{q,r} \qquad\qquad (u \in {\cal F}_{q,r}(A)).
\end{equation}

%%%%%%%%%%%%%%%%%%%%%%%%%%%%%%%%%%%%%%%%%%%%%%%
%%%%%%%%%%%%%%%%%%%%%%%%%%%%%%%%%%%%%%%%%%%%%%%
%%%%      Section      Quasi-sure lift of BM
%%%%%%%%%%%%%%%%%%%%%%%%%%%%%%%%%%%%%%%%%%%
%%%%%%%%%%%%%%%%%%%%%%%%%%%%%%%%%%%%%%%%%%%%%%%
%%%%%%%%%%%%%%%%%%%%%%%%%%%%%%%%%%%%%%%%%%%%%%%

\section{Quasi-sure existence of Brownian rough path}

In this section, we recall that Brownian motion admits a natural lift quasi-surely
via the dyadic partitions.
This fact was proved by three authors independently, Higuchi \cite{hi}, Inahama \cite{in1}, 
and Watanabe \cite{wa2}.
Among them, Higuchi's method seems best.
Higuchi's master thesis is in Japanense and probably unavailable outside Japan.
However, Section 3 of Aida's recent paper \cite{ai} is essentially the same.
So we will follow \cite{hi, ai},
in which a slightly different Besov norm is used.

From now on, we denote by $({\cal W}, {\cal H}, \mu)$ 
be the $d$-dimensional classical Wiener space.
That is, 
${\cal W} = C_0 ([0,1], {\bf R}^d)$ with the sup-norm, 
${\cal H}$ is the Cameron-Martin space, and $\mu$ is the usual $d$-dimensional Wiener measure.
For $w \in {\cal W}$ and $k \in {\bf N}$, 
$w(k) \in C_0^{1-H} ([0,1], {\bf R}^d)$ denotes the $m$th dyadic polygonal approximation 
associated with the partition $\{ l2^{-k} ~|~0 \le l \le 2^k \}$ of $[0,1]$. 
We denote by $W(k) := S_2(w(k))$ the natural lift of $w(k)$.

For $\al, m$ such that
$4m \ge 2$, $1/3 <\al <1/2$, and $\al -1/(4m) >1/3$,
we set 
\[
{\cal Z}_{\al, 4m} = \bigl\{  w \in {\cal W}  ~|~ 
\mbox{ $\{ W(k)\}_{k=1}^{\infty}$ is Cauchy in $G\Omega_{\al, 4m}^B ({\bf R}^d)$} \bigr\}.
\]
Slightly abusing the notation, we denote $\lim_{k \to \infty} W(k)$
by $W$ or $S_2(w)$ for $w \in {\cal Z}_{\al, 4m}$.
If $w \notin {\cal Z}_{\al, 4m}$, $S_2(w)$ is not defined.
(So, as a map from ${\cal W}$, $S_2$'s definition depends on $\al, m$.)
The aim of this section is to prove that, for sufficiently large $m \in {\bf N}$,
${\cal Z}_{\al, 4m}^c$ is a slim set.

Some basic properties of ${\cal Z}_{\al, 4m}$ are shown in the following lemma.
We can see from this that we may write $c W~(c \in {\bf R})$ without ambiguity.
%
%
%%%%%%%%%%%%%%  Lemma  %%%%%%%%%%%%%%%%%%%%
%
%
\begin{lm}\label{lm.Zjimei}
Let $\al, m$ and ${\cal Z}_{\al, 4m}$ be as above. Then, we have the following {\rm (i)--(iii)}.
\\
\noindent {\rm (i)}~
${\cal H} \subset {\cal Z}_{\al, 4m}$.
\\
\noindent {\rm (ii)}~
For any $c \in {\bf R}$ and $w \in {\cal Z}_{\al, 4m}$, $S_2(cw) = cS_2(w)$.
In particular, $c {\cal Z}_{\al, 4m}= {\cal Z}_{\al, 4m}$ if $c \neq 0$.
\\
\noindent {\rm (iii)}~
Assume in addition that $2m -4m\al -1 = 4m (1/2 -\al) -1 >0$.
Then, 
for any $h \in {\cal H}$ and $w \in {\cal Z}_{\al, 4m}$,  
$S_2(w +h) =  T_h (S_2(w))$.
Here, the right hand side stands for the Young translation of $S_2(w)$ by $h$.
In particular,  
$T_h ({\cal Z}_{\al, 4m} )= {\cal Z}_{\al, 4m}$
\end{lm}

\Proof
Throughout this proof,  the constant $C>0$ may change from line to line. 
{\rm (ii)} is trivial. 
Note that $\|h\|_{1/2 -H} \le \| h\|_{{\cal H}}$
and that, for any $h$, $\| h(k) \to h\|_{{\cal H}} \to 0$ as $k \to \infty$.
Then, {\rm (i)} is immediate from these.

We now prove {\rm (iii)}. 
The only non-trivial part  is boundedness of the "cross terms" of the second level path
as bilinear functionals in $X^1$ and $h$.
For $x \in C^{\al, 4m -B}_0 ([0,1], {\bf R}^d)$ and $h \in {\cal H}$,
consider 
\[
J[x,h]_{s,t} := \int_s^t  (x_u -x_s) \otimes dh_u.
\]
Then, we can easily see from (\ref{inj.HB.pr}) that
\begin{align}
|J[x,h]_{s,t}|^2 &\le  \bigl( \int_s^t  |x_u -x_s| |h^{\prime}_u | du  \bigr)^2
\le
\|h\|^2_{{\cal H}}
\int_s^t  |x_u -x_s|^2 du   
\nn\\
& \le 
C \|h\|^2_{{\cal H}} \| x\|^2_{\al, 4m-B}
\int_s^t  (u -s)^{2(\al -1/4m)  } du   
\nn\\
&\le
C \|h\|^2_{{\cal H}} \| x\|^2_{\al, 4m-B}   (t -s)^{2(\al -1/4m)  +1}
\nn
\end{align}
%
%Throughout this proof,  the constant $C>0$ may change from line to line. 
%
Then, we have
\begin{align}
\|J[x,h] \|^{2m}_{2\al, 2m -B}  
&\le
C \|h\|^{2m}_{{\cal H}}   \| x\|^{2m}_{\al, 4m-B}
\iint_{0\le s \le t \le 1} \frac{ (t -s)^{2m(\al -1/4m)  +m} }{  (t -s)^{1+4m\al}  } dsdt
\nn\\
&\le
C \|h\|^{2m}_{{\cal H}}   \| x\|^{2m}_{\al, 4m-B}
\iint_{0\le s \le t \le 1} \frac{ 1 }{  (t -s)^{1+2m\al -m +1/2}  } dsdt 
\nn\\
&
\le C \|h\|^{2m}_{{\cal H}}   \| x\|^{2m}_{\al, 4m-B}.
\nn
\end{align}
In the same way, $J[h, x]$ can be estimated.
The rest is a routine and we omit it.
\QED

%%%%%%%%%%%%%%%%%%%%%%%%%%%%%%%%%%%%%%%%%%%%%%%%%%%%%%%%%%%%%%%%%%%%%%%%%%%%%%%%
%%%   ����???����?��
%%%%%%%%%%%%%%%%%%%%%%%%%%%%%%%%%%%%%%%%%%%%%%%%%%%%%%%%%%%%%%%%%%%%%%%%%%%%%%

%%%%%%%%%%%%

Take any $\ve >0$ and fix it. We set 
\begin{align}
{\cal A}^k_{\al, 4m} = \bigl\{  w \in {\cal W}  ~|~
  \| w(k+1) -w(k) \|_{\al, 4m-B}^{4m}  \ge k^{- (4m+\ve)}
 \bigr\}, 
 \label{Asub.def}
\end{align}
Note that $\| w(k+1) -w(k) \|_{\al, 4m-B}^{4m}$ is in the  
inhomogeneous Wiener chaos $\hat{\cal C}_{4m}$ if $m \in {\bf N}$.
We can estimate the capacities of ${\cal A}^k_{\al, 4m}$ as follows.
%
%
%%%%%%%%    Lemma    %%%%%%%%%%%%%%%%%%%%%%%%%%%%%%%
%
\begin{lm}\label{lm.capA}
As before, let $m \in {\bf N}_+$, $1/3 <\al <1/2$, and $\al - (1/4m) >1/3$.
Furthermore, we assume that
$4m -8m\al -1 =8m( 1/2 - \al ) -1 >0$. 
Then, for any $1 <q <\infty$ and $r \in {\bf N}$, 
\begin{equation}
\sum_{k=1}^{\infty}
{\rm c}_{q,r} ({\cal A}^k_{\al, 4m}  )  
\le M_{q,r} 
\sum_{k=1}^{\infty} k^{4m+\ve}
\Bigl\|  \| w(k+1) -w(k) \|_{\al, 4m-B}^{4m} \Bigr\|_{q,r}  <\infty.
\label{capA.ineq}
\end{equation}
Here, $M_{q,r}>0$ is the constant given in Lemma \ref{lm.cheby}.
As a result, $\limsup_{k \to \infty} {\cal A}^k_{\al, 4m}$ is slim, 
which implies quasi-sure convergence of 
$\{ W(k)^1\}_k$. 
\end{lm}

\Proof
The left inequality in (\ref{capA.ineq}) is immediate from Lemma \ref{lm.cheby}.
First we will prove the summability for $L^2$-norm, instead of ${\bf D}_{q,r}$-norm.

We set $z(k)=  w(k+1) -w(k)$ and $Z(k)^1_{s,t} = z(k)_t - z(k)_s$ for simplicity.
It can be written explicitly as follows;
\begin{align}
z(k)_t &= 2^k \{  (t -\frac{j-1}{2^k}) \wedge (\frac{j}{2^k} -t)  \}
\nn\\
&  \qquad\qquad \times
\Bigl( 2 w (\frac{2j-1}{2^{k+1}})  -w (\frac{j-1}{2^{k}}) - w (\frac{j}{2^{k}})
\Bigr)
\qquad
\mbox{ if }   \frac{j-1}{2^{k} } \le t \le \frac{ j }{2^{k}}.
\label{zk-hyoji.eq}
\end{align}
Note that $ \| z(k) \|_{\al, 4m-B}^{4m} \in \hat{\cal C}_{4m}$.
A straight forward computation yields;
\begin{equation}\label{zkL2.ineq}
{\mathbb E} [ | Z(k)^1_{s,t} |^2 ] \le C_d \bigl( \frac{1}{2^k} \wedge (t-s) \bigr),
\qquad \qquad
 (s,t) \in \triangle
\end{equation}
for some constant $C_d >0$, which depends only on $d$.
Then, we have 
\begin{align}
\Bigl\|  \| z(k) \|_{\al, 4m-B}^{4m} \Bigr\|^2_{2}  
&\le 
{\mathbb E} \Bigl[  \Bigl( \iint_{0 \le s \le t \le 1}  
\frac{|Z(k)^1_{s,t} |^{4m}}{(t-s)^{1 +4m\al} } dsdt \Bigr)^2 \Bigr]
\nn\\
&\le
\frac 12  \iint_{0 \le s \le t \le 1}
\frac{  {\mathbb E}[ |Z(k)^1_{s,t} |^{8m} ] }{(t-s)^{2 +8 m\al} } dsdt
\nn\\
&\le
C_{m,d} 
\iint_{0 \le s \le t \le 1}
\frac{ (2^{-k})^{4m } \wedge (t-s)^{4m}  }{(t-s)^{2 +8 m\al} } dsdt.
\label{zk1.ineq}
\end{align}
Here, we used Schwarz's inequality and 
Proposition \ref{chaos.nm.pr} for $Z(k)^1_{s,t} \in {\cal C}_1$.

Now use the following well-known inequality; 
\begin{equation}\label{tknb.ineq}
\iint_{0 \le s \le t \le 1}
\frac{ \eta \wedge (t-s)^a  }{(t-s)^b} 
dsdt 
\le 
\frac{a \eta^{(a-b+1)/a} }{ (a-b+1)(b-1)}
\qquad
(0 <\eta \le 1, b>1, a>b-1).
\end{equation}
To check this formula, just change variables by $u=s, \, v=t-s$.
Then, the domain of integral becomes $\{ 0< u<1, 0<v<1, u+v<1 \}$.
The rest is easy.

Observe that,
if $4m -8m\al -1 =8m( 1/2 - \al ) -1 >0$, we can use (\ref{tknb.ineq}).
We have
\[
\Bigl\|  \| z(k) \|_{\al, 4m-B}^{4m} \Bigr\|_{2}  
\le 
C_{\al, m,d}  \Bigl( \frac{1}{ 2^{(4m -8m\al -1) /2} } \Bigr)^{k}.
\]
This proves the summability for $L^2$-norm.
As for ${\bf D}_{q,r}$-norm,
we can easily see from Proposition \ref{chaos.nm.pr} that
\begin{align}
\sum_{k=1}^{\infty} k^{4m +\ve}  \Bigl\|  \| z(k) \|_{\al, 4m-B}^{4m} \Bigr\|_{q,r} 
&\le 
M_{4m, q,r}
\sum_{k=1}^{\infty} k^{4m +\ve}  \Bigl\|  \| z(k) \|_{\al, 4m-B}^{4m} \Bigr\|_{2}  
\nn\\
&\le
C_{q,r, \al, m,d}  \sum_{k=1}^{\infty} k^{4m +\ve}  
\Bigl( \frac{1}{ 2^{(4m -8m\al -1) /2} } \Bigr)^{k}
<\infty.
\end{align}
Here, $c_{q,r, \al, m,d}$ is a positive constant which depends only on $q,r, \al, m,d$.

It is immediate from 
Lemma \ref{lm.borcan}  that $\limsup_{k \to \infty} {\cal A}^k_{\al, 4m}$ is slim.
Note that 
$w \in (\limsup_{k \to \infty} {\cal A}^k_{\al, 4m})^c 
= \cup_{N=1}^{\infty} \cap_{k=N}^{\infty} ({\cal A}^k_{\al, 4m})^c$ implies 
that
$\{ w(k)\} =\{ W(k)^1 \}$ is convergent in $(\al, m)$-Besov norm.
\QED

%%%%%%%%%%%%%%%%%%%% 2nd level %%%%%%%%%%%%%%%%%%%%%%%%%%%%%%

Next, 
let us consider the second level paths.
For ${\bf R}^d$-valued continuous paths $x$ and $y$, 
we define $J[x,y] : \triangle \to ({\bf R}^d )^{\otimes 2}$ by
\[
J[x,y]_{s,t} = \int_s^t  (x_u-x_s) \otimes dy_u,
\]
whenever the integral on the right hand side can be defined.
Note that the following equality holds; 
\begin{align}
J[x,x]_{s,t} - J[y,y]_{s,t} 
&=
J[x-y,x-y]_{s,t} + J[x-y,y]_{s,t} 
\nn\\
& \qquad
- J[x-y,y]_{s,t}^* 
+ Y^1_{s,t} \otimes ( X^1_{s,t} -  Y^1_{s,t}).
\label{Jalg.eq}
\end{align}
Here, $*$ stands for the linear isometry of $({\bf R}^d )^{\otimes 2}$ defined by 
$(\xi \otimes \eta)^* = \eta \otimes \xi$ for all $\xi, \eta \in {\bf R}^d$.
Indeed, we can easily see that
\[
 J[x-y,x-y]_{s,t} -J[x,x]_{s,t} + J[y,y]_{s,t} 
=
-J[x-y,y]_{s,t} - J[y, x-y]_{s,t}
\]
and, from integration by parts, that
\[
J[y, x-y]_{s,t} = Y^1_{s,t} \otimes ( X^1_{s,t} -  Y^1_{s,t})
- \int_s^t   dy_u \otimes ( X^1_{s,u} -  Y^1_{s,u}).
\]
Thus, we have verified (\ref{Jalg.eq}).

With (\ref{Jalg.eq}) in hand, 
we naturally consider the following subsets in ${\cal W}$. 
\begin{align}
{\cal B}(1)^k_{2\al, 2m} 
&= \bigl\{  w \in {\cal W}  ~|~
  \| J[z(k), z(k)] \|_{2\al, 2m-B}^{2m}  \ge k^{- (2m+\ve)}
 \bigr\}, 
 \label{B1sub.def}
 \\
 {\cal B}(2)^k_{2\al, 2m} 
 &= \bigl\{  w \in {\cal W}  ~|~
  \| J[z(k), w(k)] \|_{2\al, 2m-B}^{2m}  \ge k^{- (2m+\ve)}
 \bigr\}, 
 \label{B2sub.def}
 \\
 {\cal B}(3)^k_{2\al, 2m} 
 &= 
 \bigl\{  w \in {\cal W}  ~|~
  \| W(k)^1 \otimes  Z(k)^1 \|_{2\al, 2m-B}^{2m}  \ge k^{- (2m+\ve)}
 \bigr\},
 \label{B3sub.def}
  \end{align}
where we set $z(k) =w(k+1) -w(k)$ as before.

%%%%%%%%%%    Lemma    %%%%%%%%%%%%%%%%%%%%%%%%%%
%
%
\begin{lm}\label{lm.capB}
As before, let $m \in {\bf N}_+$, $1/3 <\al <1/2$, and $\al -(1/4m) >1/3$.
Furthermore, we assume that $4m -8m\al -1 =8m( 1/2 - \al ) -1 >0$.
Then, for any $1 <q <\infty$ and $r \in {\bf N}$, 
\begin{align}
\sum_{k=1}^{\infty}
{\rm c}_{q,r} ({\cal B}(1)^k_{2\al, 2m}  )  
&\le M_{q,r} 
\sum_{k=1}^{\infty} k^{2m+\ve}
\Bigl\|  \| J[z(k), z(k)] \|_{2\al, 2m-B}^{2m} \Bigr\|_{q,r}  <\infty,
\label{capB1.ineq}
\\
\sum_{k=1}^{\infty}
{\rm c}_{q,r} ({\cal B}(2)^k_{2\al, 2m}  )  
&\le M_{q,r} 
\sum_{k=1}^{\infty} k^{2m+\ve}
\Bigl\|  \| J[z(k), w(k)] \|_{2\al, 2m-B}^{2m} \Bigr\|_{q,r}  <\infty,
\label{capB2.ineq}
\\
\sum_{k=1}^{\infty}
{\rm c}_{q,r} ({\cal B}(3)^k_{2\al,2m}  )  
&\le M_{q,r} 
\sum_{k=1}^{\infty} k^{2m+\ve}
\Bigl\| \| W(k)^1 \otimes  Z(k)^1   \|_{2\al, 2m-B}^{2m} \Bigr\|_{q,r}  <\infty.
\label{capB3.ineq}
\end{align}
Here, $M_{q,r}>0$ is the constant given in Lemma \ref{lm.cheby}.
As a result, $\limsup_{k \to \infty} {\cal B}(i)^k_{2\al, 2m}$ is slim for all $i=1,2,3$, 
which implies quasi-sure convergence of  $\{ W(k)^2\}_k$. 
\end{lm}

\Proof
Substituting $x=w(k+1)$ and $y=w(k)$ in (\ref{Jalg.eq}), 
we  have
\[
W(k+1)^2 - W(k)^2
=
J[z(k), z(k)] + J[z(k), w(k)] -  J[z(k), w(k)]^* + W(k)^1 \otimes  Z(k)^1.
\]
If $w \in \cap_{1 \le i \le 3} (\limsup_{k \to \infty} {\cal B}(i)^k_{\al, m})^c$,
then $\{ W(k)^2 \}_{k=1}^{\infty}$ 
is a Cauchy sequence in $(2\al, 2m)$-Besov norm.
The left inequalities in (\ref{capB1.ineq})--(\ref{capB3.ineq})
are obvious from Lemma \ref{lm.cheby} of Chebyshev type.
By Lemma \ref{lm.borcan} of Borel-Cantelli type, 
summability in (\ref{capB1.ineq})--(\ref{capB3.ineq}) implies that 
$\limsup_{k \to \infty} {\cal B}(i)^k_{\al, m}$ are slim for all $i$.

Hence, it is sufficient to show that the sums on the right hand sides of 
(\ref{capB1.ineq})--(\ref{capB3.ineq}) converge.
First, we consider (\ref{capB1.ineq}).
It is not difficult to see that, for a positive constant $C_d$, 
\begin{equation}\label{l2jzz.eq}
{\mathbb E}[| J[z(k), z(k)]_{s,t}|^2] 
 \le C_d 
 \bigl( \frac{1}{ 2^{2k} }  \wedge   (t-s)^2 \bigr),
 \qquad 
 \quad
 (s,t) \in \triangle
\end{equation}
A rough sketch of proof of this estimate is as follows.
Suppose that $(j-1)/2^{k} \le s \le j /2^{k}$
and $(l-1)/2^{k} \le t \le l /2^{k}$ with $j \le l$.
If $j=l$, then we have from (\ref{zk-hyoji.eq}) that
\begin{align}
J[z(k) ,z(k)]_{s,t} &= \frac{ 2^{2k}}{2} 
\bigl\{  (t -\frac{j-1}{2^k}) \wedge (\frac{j}{2^k} -t)  -  (s -\frac{j-1}{2^k}) \wedge (\frac{j}{2^k} -s)  \bigr\}^2
\nn\\
&  \qquad\qquad \times
\Bigl( 2 w (\frac{2j-1}{2^{k+1}})  -w (\frac{j-1}{2^{k}}) - w (\frac{j}{2^{k}})
\Bigr)^{\otimes 2}.
%\label{jzz-hyoji.eq}
\nn
\end{align}
In this case, (\ref{l2jzz.eq}) is easy.
If $j <l$, we use Chen's identity to obtain
\begin{align}
J[z(k) ,z(k)]_{s,t} &=J[z(k) ,z(k)]_{s, j /2^{k}}+ J[z(k) ,z(k)]_{(l-1)/2^{k} ,t}
+ Z(k)^1_{s, j /2^{k}}  \otimes Z(k)^1_{ (l-1)/2^{k} , t}.
\nn
\end{align}
In the same way as above, we can estimate the right hand side. Thus, we have shown (\ref{l2jzz.eq}).

With (\ref{l2jzz.eq}) in hand, we can calculate in the same way as in (\ref{zk1.ineq})
to obtain
\[
\Bigl\|  \| J[z(k),z(k)] \|_{2\al, 2m-B}^{4m} \Bigr\|_{2}  
\le 
C_{\al, m,d}  \Bigl( \frac{1}{ 2^{(4m -8m\al -1) /2} } \Bigr)^{k}.
\]
From this and Proposition \ref{chaos.nm.pr}, we can easily see that,
if $4m -8m\al -1 >0$,
the sum in (\ref{capB1.ineq}) converges.

Next, we consider (\ref{capB3.ineq}).
In this case we have
\begin{equation}\label{l2zw.eq}
{\mathbb E}[|   W(k)^1_{s,t} \otimes  Z(k)^1_{s,t} |^2] 
 \le C_d 
 \bigl( \frac{1}{ 2^{k} }  \wedge   (t-s)^2  \bigr),
 \qquad 
 \quad
 (s,t) \in \triangle
\end{equation}
for some positive constant $C_d$.
From (\ref{zkL2.ineq}) and an estimate that ${\mathbb E}[|   W(k)^1_{s,t}|^2] \le d(t-s)$ 
for all $k$,  we can easily show (\ref{l2zw.eq}).

From this and 
Proposition \ref{chaos.nm.pr} for $W(k)^1_{s,t} \otimes Z(k)^1_{s,t} \in \hat{\cal C}_2$,  we have 
\begin{align}
\Bigl\|  \| W(k)^1 \otimes  Z(k) \|_{2\al, 2m-B}^{2m} \Bigr\|^2_{2}  
&\le 
{\mathbb E} \Bigl[  \Bigl( \iint_{0 \le s \le t \le 1}  
\frac{| W(k)^1_{s,t} \otimes Z(k)^1_{s,t} |^{2m}}{(t-s)^{1 +4m\al} } dsdt \Bigr)^2 \Bigr]
\nn\\
&\le
\frac 12  \iint_{0 \le s \le t \le 1}
\frac{  {\mathbb E}[ |  W(k)^1_{s,t} \otimes Z(k)^1_{s,t} |^{4m} ] }{(t-s)^{2 +8 m\al} } dsdt
\nn\\
&\le
C_{m,d} 
\iint_{0 \le s \le t \le 1}
\frac{ (2^{-k})^{2m } \wedge (t-s)^{4m}  }{(t-s)^{2 +8 m\al} } dsdt
\nn\\
&
\le
C_{\al, m,d}   \Bigl(  \frac{1}{ 2^{ ( 4m -8m\al -1) /2} } \Bigr)^k.
\label{w1z1.ineq}
\end{align}
Here, we used (\ref{tknb.ineq}) for the last inequality.

Applying Proposition \ref{chaos.nm.pr} to
$\| W(k)^1 \otimes  Z(k) \|_{2\al, 2m-B}^{2m} \in \hat{\cal C}_{4m}$, 
we obtain  from (\ref{w1z1.ineq}) that 
\[
\Bigl\|  \| W(k)^1 \otimes  Z(k) \|_{2\al, 2m-B}^{2m} \Bigr\|_{q,r}
\le 
C_{q,r, \al, m,d}   \Bigl(  \frac{1}{ 2^{ ( 4m -8m\al -1) /4} } \Bigr)^k
\]
for some positive constant $C_{q,r, \al, m,d}$ independent of $k$.
Hence, if $4m -8m\al -1 >0$,  the right hand side of (\ref{capB3.ineq}) is summable.

Finally, we consider (\ref{capB2.ineq}).
In this case we have
\begin{equation}\label{l2Jzw.eq}
{\mathbb E}[|   J[ z(k),  w(k) ]_{s,t} |^2] 
 \le C_d 
 \bigl( \frac{1}{ 2^{k} }  \wedge   (t-s)^2  \bigr),
 \qquad 
 \quad
 (s,t) \in \triangle
\end{equation}
for some positive constant $C_d$.
The right hand side of (\ref{l2Jzw.eq}) is the same as  in (\ref{l2zw.eq}).
Hence,  (\ref{l2Jzw.eq})  implies summability of  (\ref{capB2.ineq}).

Now we give a sketch of proof of  (\ref{l2Jzw.eq}).
If $(j-1)/2^{k} \le s \le t \le j /2^{k}$ for some $j$, 
then we can easily verify (\ref{l2Jzw.eq}) by using (\ref{zkL2.ineq}).
If $(j-1)/2^{k} \le s \le j /2^{k}$
and $(l-1)/2^{k} \le t \le l /2^{k}$ with $j < l$, then 
we can verify (\ref{l2Jzw.eq}) by using Chen's identity as follows; 
\begin{align}
J[z(k) ,w(k)]_{s,t} 
&=
J[z(k) ,w(k)]_{s, j /2^{k}}
+ 
\sum_{ j+1 \le i \le l-1} J[z(k) ,w(k)]_{ (i-1) /2^{k}, i /2^{k}}
\nn\\
& \qquad
+
J[z(k) ,w(k)]_{(l-1)/2^{k} ,t}
+ Z(k)^1_{s, j /2^{k}}  \otimes W(k)^1_{ j/2^{k} , t}.
\nn
\end{align}
From this we can see (\ref{l2Jzw.eq}).
Thus, we have shown the lemma.
 \QED

As a corollary, we have quasi-continuity of Brownian rough path,
which was first given by Aida \cite{ai} in a slightly stronger form.

\begin{co}\label{co.qcontW}
%As before, let $m \in {\bf N}_+$, $1/3 <\al <1/2$, and $\al - (1/4m) >1/3$.
%Furthermore, we assume that
%$4m -8m\al -1 =8m( 1/2 - \al ) -1 >0$. 
%
Let $\al$ and $m$ be as in Lemmas \ref{lm.capA} and \ref{lm.capB}.
Then, 
%${\cal Z}_{\al, 4m}^c$ is slim and 
the mapping 
$w \mapsto W= S_2(w)$ is $\infty$-quasi continuous.
\end{co}

\Proof
%Everything except quasi-continuity was proved.
We will show $S_2: {\cal W} \to G\Omega^B_{\al, 4m} ({\bf R}^d)$ is quasi-continuous.
Let $n \in {\bf N}_+$ be arbitrary.
Then, we can find $N=N(n)$ such that 
\[
\sum_{k=N}^{\infty}
{\rm c}_{n,n} ({\cal A}^k_{\al, 4m}  ) < \frac{1}{4n},
\quad
\max_{1 \le i \le 3}
\sum_{k=N}^{\infty}
{\rm c}_{n,n} ({\cal B}(i)^k_{2\al, 2m}  ) <\frac{1}{4n}.
\]
On $(\cup_{k=N}^{\infty} {\cal A}^k_{\al, 4m})^c =\cap_{k=N}^{\infty} ({\cal A}^k_{\al, 4m})^c$, 
$\{ W^1 (k)\}_k$ converges uniformly as $k \to \infty$.
For each $k$, $w \mapsto W^1 (k)$ is clearly continuous.
Hence, outside $\cup_{k=N}^{\infty} {\cal A}^k_{\al, 4m}$,
$W^1 = S_2(w)^1$ is continuous in $w$.
Similarly, 
outside $\cup_{i=1}^3 \cup_{k=N}^{\infty} {\cal B}(i)^k_{\al, 4m}$,
$W^2 = S_2(w)^2$ is continuous in $w$.

Therefore, the set
$\cup_{k=N}^{\infty} {\cal A}^k_{\al, 4m} \cup 
[\cup_{i=1}^3 \cup_{k=N}^{\infty} {\cal B}(i)^k_{\al, 4m}]$ 
has $(n,n)$-capacity
smaller than $1/n$, due to subadditivity of the capacity,
and its complement is a subset of ${\cal Z}_{\al, 4m}$.
Moreover, on its complement, $S_2$ is continuous.
This proves the corollary.
\QED

\begin{re}\label{re.am}
If we assume (\ref{eq.amam})
%
%\begin{equation}\label{eq.amam}
%m \in {\bf N}_+, \quad  \frac13 <\al < \frac12, \quad \al - \frac{1}{4m} > \frac13, \quad
%\mbox{ and } \quad 4m -8m\al   >2,
%\end{equation}
%
then all the assumptions for $\al, m$ in all the results in this section are satisfied.
The condition (\ref{eq.amam}) might not be best possible, but is sufficient for our purpose.
Because what is really needed is as follows; 
for any $\al \in (1/3, 1/2)$, sufficiently large $m \in {\bf N}_+$ satisfies 
the assumptions of these results.
\end{re}

%%%%%%%%%%%%%%%%%%%%%%%%%%%%%%%%%%%%%%%%%%%%%%%%%%%%%%%%%%%%%%%%%%%%%%%%%%%%%%%%%%%%%%%%%%%%%%
%%%%%%%%%%%%%%%%%%%%%%%%%%%%%%%%%%%%%%%%%%%%%%%%%%%%%%%%%%%%%%%%%%%%%%%%%%%%%%%%%%%%%%
%%%%%%%%%%%%%%%%%%%%%%%%%%%%%%%%%%%%%%%%%%%%%%%%%%%%%%%%%%%%%%%%%%%%%%%%%%%%%%%%%%%%
%%  Section    SDE  
%%%%%%%%%%%%%%%%%%%%%%%%%%%%%%%%%%%%%%%%%%%%%%%%%%%%%%%%%%%%%%%%%%%%%%%%%%%%%%%%%%%%%%
%%%%%%%%%%%%%%%%%%%%%%%%%%%%%%%%%%%%%%%%%%%%%%%%%%%%%%%%%%%%%%%%%%%%%%%%%%%%%%%%%%%%%%
%%%%%%%%%%%%%%%%%%%%%%%%%%%%%%%%%%%%%%%%%%%%%%%%%%%%%%%%%%%%%%%%%%%%%%%%%%%%%%%%%

\section{Quasi-sure analysis of SDE via rough path}
\label{se.six}

We study SDE (\ref{sde_set.eq}),
which involves a small positive parameter $\ve \in (0,1]$.
In this section  $\ve$ is fixed, however. 
Hence, readers who like it simple may set $\ve =1$.
%
%

%Let $\sigma: {\bf R}^n \to {\rm Mat}(n,d)$ be $C_b^3$, that is, $\sigma$ is 
%a bounded and $C^3$ function whose derivatives are all bounded.
%A $C_b^3$ function
%$b: [0,1] \times {\bf R}^n \to {\bf R}^n$ is an $\ve$-dependent vector field.
%
%
%The $i$th column vector of $\sigma$ is denoted by $V_i : {\bf R}^n \to {\bf R}^n$ ($1 \le i \le d$).
%When $V_i$'s are used,  $b(\ve, a)$ is denoted by $V_0 (\ve, a)$.

%We consider the following ${\bf R}^n$-valued SDE of Stratonovich-type;
%\begin{equation}\label{sde.def}
%dy^{\ve}_t  = \ve \sigma( y^{\ve}_t) \circ  dw_t + b(\ve,  y^{\ve}_t)dt,
%\qquad\qquad
% y^{\ve}_0 = a \in {\bf R}^n.
%\end{equation}
%
%
%Here, $(w_t)$ is the canonical realization of $d$-dimensional Brownian motion.
%Equivalently, (\ref{sde.def}) can also be written as follows;
%\begin{equation}\label{sde2.def}
%dy^{\ve}_t  = \ve  \sum_{i=1}^d  V_i( y^{\ve}_t) \circ  dw_t^i +  V_0 (\ve,  y^{\ve}_t)dt, 
%\qquad\qquad  y^{\ve}_0 = a.
%\end{equation}
%

Next we introduce a rough differential equation (RDE)
associated with the vector fields $V_i$'s
and express $y^{\ve}$ in terms of  the Lyons-It\^o map.
For a while we assume $V_i$'s are only $C_b^3$.
Set $\lambda \in C_0^{1 - H} ([0,1], {\bf R})$ by $\lambda_t =t$.
For $1/3 < \al <1/2$,
we denote by $\tau_{\lambda} : G\Omega^H_{\al} ({\bf R}^d) \to G\Omega^H_{\al} ({\bf R}^{d+1})$ 
the Young pairing formally given by $X \mapsto (X, \lambda)$.
We denote by $\Phi_{\ve} : G\Omega^H_{\al} ({\bf R}^{d+1}) \to G\Omega^H_{\al} ({\bf R}^{n})$
the Lyons-It\^o map associated with 
the coefficients  $[V_1, \ldots , V_d,V_0(\ve, \,\cdot\,)]$
and the starting point $a$,
which is continuous 
by Lyons' continuity theorem. 
If $X =S_2(x)$ for a $1$-H\"older continuous path, 
then $a+\Phi_{\ve} (\tau_{\lambda} (X) )^1$ 
solves the following ODE in Riemann-Stieltjes sense; 
\[
dz_t  =  \sum_{i=1}^d V_i ( z_t)   dx_t^i + V_0 (\ve, z_t)dt,
\qquad\qquad
 z_0 = a.
\]
Hence, by Wong-Zakai's approximation theorem, 
$y^{\ve} =a+ \Phi_{\ve} (\tau_{\lambda} (\ve W) )^1$ for a.a.$w$.
Thus, the solution of the scaled SDE is obtained via rough path theory.
Note also that this shows $y^{\ve}$ is a $C^{\al -H} ([0,1], {\bf R}^n)$-valued 
Wiener functional.

\begin{pr}\label{pr.MalNua}
Assume that $V_1, \ldots , V_d,V_0(\ve, \,\cdot\,)$ are $C_b^3$
and let
$(y^{\ve}_t)_{0 \le t \le 1}$ be the solution of SDE (\ref{sde_set.eq}).
Then, for any $1/3 < \al <1/2$, the mapping
\[
{\cal W} \ni w \mapsto y^{\ve} \in C^{\al -H} ([0,1], {\bf R}^n)
\]
admits an $\infty$-quasi-redefinition, which is explicitly given by
 $a+\Phi_{\ve} (\tau_{\lambda} (\ve W) )^1$.
\end{pr}

\Proof
For given $\al$, we can find $(\al', m)$ as in (\ref{eq.amam}) such that $\al <\al' $
and $\al' - (4m)^{-1} \ge \al$.
By Proposition \ref{inj.HB.pr},  $G\Omega^B_{\al', 4m} ({\bf R}^d)$ is continuously 
embedded in $G\Omega^H_{\al} ({\bf R}^d)$.
To keep  notations simple, we will not write this injection explicitly.

For all $w \in {\cal W}$ and $k \in {\bf N}_+$,  
we have $y(k)^{\ve} = a+ \Phi_{\ve} (\tau_{\lambda} (\ve W(k) ) )^1$,
where $y(k)^{\ve}$ is the unique  solution of 
\[
d  y(k)^{\ve}_t  =  \sum_{i=1}^d V_i ( y(k)^{\ve}_t)  \ve  d w(k)_t + V_0 (\ve,  y(k)^{\ve}_t)dt,
\qquad\qquad
 y(k)^{\ve}_0 = a.
\]
Since we assume $C_b^3$, 
we may use Wong-Zakai's theorem. 
It says $\sup_{0 \le t \le 1} |  y(k)^{\ve}_t -  y^{\ve}_t | \to 0$ almost surely
as $k \to \infty$.
On the other hand, by Lyons' continuity theorem, 
Proposition \ref{inj.HB.pr}, and Corollary \ref{co.qcontW}, 
$\Phi_{\ve} (\tau_{\lambda} (\ve W(k) ) )^1$ converges to 
$\Phi_{\ve} (\tau_{\lambda} (\ve W ) )^1$ in $\al$-H\"older norm
not just almost surely, but also quasi-surely.
%
%(Note that there exist $\al' \in (\al, 1/2)$ and $m \in {\bf N}_+$ 
%which satisfy the assumption of  Corollary \ref{co.qcontW} and $\al' - (4m)^{-1} \ge \al$.)
%
%
Hence, we have $y^{\ve} = a+\Phi_{\ve} (\tau_{\lambda} (\ve W ) )^1$ a.s.
Again by Lyons' continuity theorem, Proposition \ref{inj.HB.pr}, and Corollary \ref{co.qcontW}, 
the right hand side is $\infty$-quasi-continuous in $w$.
This implies that $a+\Phi_{\ve} (\tau_{\lambda} (\ve W ) )^1$ is an $\infty$-quasi-redefinition of  $y^{\ve}$.
\QED

\begin{re}
J. Ren \cite{ren1} studied the quasi-sure refinement of Wong-Zakai's approximation theorem.
The proof of the above theorem is  a new proof of
Ren's result.
In fact, it is an improvement since we only assumed $C_b^3$, not $C_b^{\infty}$ for the coefficients.
\end{re}

\begin{re}
The above proposition can  be regarded as a new proof 
of Malliavin-Nualart's result in \cite{mn}. 
Here, we give a quick remark about it.
Assume that the coefficients $V_i$'s are $C_b^{\infty}$.
Then, for each $t$, $y^{\ve}_t \in {\bf D}_{\infty} ({\bf R}^n)$ 
and therefore it admits an $\infty$-quasi-redefinition by a well-known general theory.
J. Ren \cite{ren1}  proved that there exists a nice modification $\hat{y}^{\ve}$ of $y^{\ve}$ 
such that, for each $t$, $\hat{y}^{\ve}_t$ is an $\infty$-quasi-redefinition of $y^{\ve}_t$.
(This process $\hat{y}^{\ve}$ is the limit of Ren's refinement of Wong-Zakai's theorem 
in \cite{ren1}.)

But, it is more difficult to prove 
the path space-valued functional  $y^{\ve}$ admits an $\infty$-quasi-redefinition.
Malliavin and Nualart \cite{mn} proved that $w \mapsto y^{\ve}$ 
admits an $\infty$-quasi-redefinition.
(It might be possible to prove this
by extending Ren's method, though it is not explicitly written in \cite{ren1}.)
More precisely, they proved that Wiener functional 
\begin{equation}\label{2para.eq}
{\cal W} \ni w \mapsto \Bigl[  [0,1] \times {\bf R}^n  \ni  (t, a)  \mapsto   y^{\ve}_t = y^{\ve}(t, a) \Bigr]
\end{equation}
which takes values in a space of two-parameter functions, admits an $\infty$-quasi-redefinition. 
But, this is not the point here.
%(It might be possible to prove that $y^{\ve}$ an $\infty$-quasi-redefinition
%by extending Ren's method in \cite{ren1})

Their method is as follows.  
Firstly, they introduce a UMD Banach space with certain nice properties 
and showed  $y^{\ve}$ takes values in it. 
Secondly, they extended Malliavin calculus so that it applies to Wiener functionals 
which takes values in the UMD Banach space.
Thirdly, they regarded SDE (\ref{sde_set.eq}) as an equation that takes values 
in the UMD Banach space, took ${\cal H}$-derivative of it to prove that $y^{\ve}$
is a UMD Banach space-valued ${\bf D}_{\infty}$-Wiener functional. 
Then, the same argument as in the case of  real-valued ${\bf D}_{\infty}$-Wiener functionals
applies to $y^{\ve}$, too, and existence of $\infty$-redefinition is shown.

Note that, in their method, smoothness of the coefficients needs to be assumed.
In our Proposition \ref{pr.MalNua}, however, we only assumed $C^3_b$, 
which is probably astonishing if we do not know rough path theory.
In this sense, this is an improvement of  Malliavin-Nualart's result.

Since Lyons-It\^o map is continuous in the initial value $a$, too, 
it is clearly possible to extend Proposition \ref{pr.MalNua} to the case of 
(\ref{2para.eq}), by introducing certain H\"older-type norm on the two-parameter function space.
It is not clear, however,  that  we have really improved the results in \cite{mn} or not,
mainly because we do not know how strong  the UMD Banach norm is, compared to 
the H\"older-type norm.
(Practically, such slight difference of norms is not important, however.)
\end{re}

%%%%%%%%%%%%%%%%%%%%%%%%%%%%%%%%%%%%%%%%%%%%%%%%%%%%%%
%%%%%%%%%%%%%%%%%%%%%%%%%%%%%%%%%%%%%%%%%%%%%%%%%%%%%%%%%%%%%%%%%%%%
%%%%%%%%%%%%%%%%%%%%%%%%%%%%%%%%%%%%%%%%%%%%%%%%%%%%%%%%%%%%%%%%%

Now, we discuss pinned diffusion measures. 
By using the arguments we developed so far,
we can obtain pinned diffusion processes via rough path theory. 
This result was first proved in Inahama \cite{in1}, but the argument here is much more 
straight forward, due to quasi-continuity of $W$ (Corollary \ref{co.qcontW}).

From now on we assume  $V_i ~(0 \le i \le d)$
are $C_b^{\infty}$.
We assume further that the vector fields $\{V_1, \ldots, V_d ; V_0^{\ve}\}$
satisfy hypoellipticity condition at any $a$.
Here, $V_0^{\ve} =V_0(\ve, \,\cdot\,)$.
Then, the solution $y^{\ve}_t = y^{\ve} (t,a)$ 
of SDE (\ref{sde_set.eq}) has a density $p_t^{\ve} (a, a')$ 
for all $t>0$, $\ve>0$, and $a \in {\bf R}^n$,
that is,  for any Borel set $A \subset {\bf R}^n$, 
$
{\mathbb P} ( y^{\ve}(t,a)  \in A) = \int_A   p_t^{\ve} (a, a') da'. $
By Watanabe's theory, $p_t^{\ve} (a, a')= {\mathbb E} [ \delta_{a'} ( y^{\ve} (t,a) ) ]$.
(See Section 5.9 --5.10, \cite{iwbk}).
In this case $p^{\ve}_1 (a,a') >0$ for any $a, a'$.

%Take any $a, a'$ such that $p^{\ve}_1 (a,a') >0$.
%

%Let ${\mathbb Q}^{\ve}_{a,a'}$ be the pinned diffusion measure from $a$ to $a'$
%associated with the second order differential operator
%${\cal L}^{\ve} = (\ve^2 /2) \sum_{i=1}^d V_i^2  + V_0^{\ve}.$
%
%
%

Under this condition, 
${\mathbb Q}^{\ve}_{a,a'}$ exists as a probability measure on 
$C_{a,a'} ([0,1], {\bf R}^n) :=\{ x \in  C([0,1], {\bf R}^n) ~|~ x_0 =a, x_1=a' \}$
characterized by the following equation;
for any $l \in {\bf N}_+$, any partition $\{ 0 =t_0 <t_1 <\cdots < t_l <t_{l+1} =1\}$ of $[0,1]$,
 and any $f \in {\cal S} ({\bf R}^{nl})$,
\[
\int f (x_{t_1}, \ldots, x_{t_l})  {\mathbb Q}^{\ve}_{a,a'}(dx)
=
p_1^{\ve} (a,a')^{-1}
\int_{ ({\bf R}^n)^l }   f (a_1, \ldots, a_l) 
\prod_{j=1}^{l+1}     p_{t_j -t_{j-1}}^{\ve} (a_{j-1}, a_j)
 \prod_{j=1}^l da_j.
\]
Here, $a_0 =a$ and $a_{l+1} =a'$ by convention.
Actually, by Proposition \ref{pr.rep_pin} below,
${\mathbb Q}^{\ve}_{a,a'}$ turns out to be a probability measure on $C^{\al -H}_{a,a'} ([0,1], {\bf R}^n)$
for any $1/3 <\al <1/2$.

Let $\theta^{\ve}_{a,a'}$ be a Borel finite measure on ${\cal W}$ that corresponds to 
$\delta_{a'} ( y^{\ve} (1,a) )$ via Sugita's theorem (Proposition \ref{pr.sugita}).
We denote by $\hat\theta^{\ve}_{a,a'}$ its normalization.
Let $\al, m$ be as in (\ref{eq.amam}).
We set $\mu^{\ve}_{a,a'} = (\ve S_2)_* \theta^{\ve}_{a,a'}$ 
and $\hat\mu^{\ve}_{a,a'} = (\ve S_2)_* \hat\theta^{\ve}_{a,a'}$.
Since ${\cal Z}_{\al, 4m}^c$ is slim, the lift map $\ve S_2$ is well-defined 
with respect to these measures.
Thus, we obtained measures on $G\Omega^B_{\al,4m} ({\bf R}^n)$.
It is easy to see that the support of 
$\hat\mu^{\ve}_{a,a'}$  is in fact contained in the closed subset 
$\{ X ~|~  a+ \Phi_{\ve} (\tau_{\lambda} (X) ) =a' \}$.

The following proposition states that 
the pinned diffusion measure can be obtained as a continuous image 
of Lyons-It\^o map.
Note that everywhere hypoellipticity was assumed to garantee existence of ${\mathbb Q}^{\ve}_{a,a'}$, 
but none of {\bf (H0)}--{\bf (H2)} was.

%This has already been shown in \cite{in1},
%but a proof given below is much more straight forward.

%%%%

\begin{pr}\label{pr.rep_pin}
Keep all the notations and assumptions as above.
Then, the push-forward measure 
$(a+\Phi_{\ve}^1 \circ \tau_{\lambda})_* \hat\mu^{\ve}_{a,a'}$
is the pinned diffusion measure ${\mathbb Q}^{\ve}_{a,a'}$. 
Here, $\Phi_{\ve}^1$ denotes the first level of the It\^o map $\Phi_{\ve}$.
\end{pr}

\Proof
For simplicity,  we write $\tilde{y}^{\ve}  = a+ \Phi_{\ve} (\tau_{\lambda} (\ve S_2(w) ) )^1$
as in the proof of Proposition \ref{pr.MalNua}.
By Proposition \ref{pr.MalNua}, this is an $\infty$-quasi-redefinition of $y^{\ve}$.
It is sufficient to prove that 
$(a+ \Phi_{\ve}^1  \circ \tau_{\lambda} \circ \ve S_2 )_* \theta^{\ve}_{a,a'}
 =p^{\ve}_1 (a,a') \cdot  {\mathbb Q}^{\ve}_{a,a'}$.
This fact is well-known in quasi-sure analysis, but we give a proof here for the reader's
convenience.

%%%

Let $\chi: {\bf R}^n \to {\bf R}$ be $C^{\infty}$ with compact support 
such that $\chi \ge 0$ and $\int_{ {\bf R}^n  } \chi (\xi) d\xi=1$.
For $k \in {\bf N}_+$, set $\psi_k (\xi) = k^n \chi (k (\xi -a'))$.
Then, for any function $g: {\bf R}^n \to {\bf R}$ which is continuous around $a'$,
we have $\lim_{k \to \infty} \int_{ {\bf R}^n  } \psi_k (\xi)  g(\xi) d\xi = g(a')$.
In particular, $\lim_{k \to \infty} \psi_k = \delta_{a'}$
 in ${\cal S}^{\prime}({\bf R}^n)$.

Take any  partition  $\{ 0 =t_0 <t_1 <\cdots < t_l <t_{l+1} =1\}$ of $[0,1]$
 and  any $f \in {\cal S} ({\bf R}^{nl})$.
By Proposition \ref{pr.sugita}, we have 
\begin{align}
\int_{{\cal W}} f (\tilde{y}^{\ve}_{t_1}, \ldots,  \tilde{y}^{\ve}_{t_l})  \theta^{\ve}_{a,a'}(dw)
&=
{\mathbb E} \bigl[  
f (y^{\ve}_{t_1}, \ldots,  y^{\ve}_{t_l})  \cdot \delta_{a'} ( y^{\ve}_t )
\bigr]
\nn\\
&
= \lim_{k \to \infty} {\mathbb E} \bigl[  
f (y^{\ve}_{t_1}, \ldots,  y^{\ve}_{t_l})  \psi_k ( y^{\ve}_t )
\bigr]
\nn\\
&
= \lim_{k \to \infty} 
\int_{ ({\bf R}^n)^{l+1} }   f (a_1, \ldots, a_l) \psi_k (a_{l+1})
\prod_{j=1}^{l+1}     p_{ t_j -t_{j-1} }^{\ve} (a_{j-1}, a_j)
 \prod_{j=1}^{l+1} da_j.
\nn\\
&
= 
\int_{ ({\bf R}^n)^{l} }   f (a_1, \ldots, a_l) 
\prod_{j=1}^{l}     p_{ t_j -t_{j-1} }^{\ve} (a_{j-1}, a_j)   p_{ t_{l+1} -t_{l} }^{\ve} (a_{l}, a')
 \prod_{j=1}^{l} da_j.
\nn
\end{align}
where $a_0 =a$.
This completes the proof.
\QED

%\newpage

%%%%%%%%%%%%%%%%%%%%%%%%%%%%%%%%%%%%%%%%%%%%%%%%%%%%%%%%%%%%%%%%%%%%%%%%%%%%%%%%%%%%%%%%%%%%%%
%%%%%%%%%%%%%%%%%%%%%%%%%%%%%%%%%%%%%%%%%%%%%%%%%%%%%%%%%%%%%%%%%%%%%%%%%%%%%%%%%%%%%%
%%%%%%%%%%%%%%%%%%%%%%%%%%%%%%%%%%%%%%%%%%%%%%%%%%%%%%%%%%%%%%%%%%%%%%%%%%%%%%%%%%%%
%%  Section    LDP  
%%%%%%%%%%%%%%%%%%%%%%%%%%%%%%%%%%%%%%%%%%%%%%%%%%%%%%%%%%%%%%%%%%%%%%%%%%%%%%%%%%%%%%
%%%%%%%%%%%%%%%%%%%%%%%%%%%%%%%%%%%%%%%%%%%%%%%%%%%%%%%%%%%%%%%%%%%%%%%%%%%%%%%%%%%%%%
%%%%%%%%%%%%%%%%%%%%%%%%%%%%%%%%%%%%%%%%%%%%%%%%%%%%%%%%%%%%%%%%%%%%%%%%%%%%%%%%%

\section{Proofs of large deviations}
\label{se.7}

%
%
%
%Under this assumption, the heat kernel $p^{\ve}_t (a,a')$ exists and positive 
%for all $a, a' \in {\bf R}^n$, $t>0$ and $\ve>0$.
%Hence, the results in the previous section are available.

The aim of this section is to prove our main results.
The most important one is Theorem \ref{tm.ldp.sgt} {\rm (i)},
which will be shown in Subsections \ref{se.low} and  \ref{se.up}.

%that 
%the family $\{ \mu^{\ve}_{a,a'}\}_{\ve >0}$ of finite measures on $G\Omega^B_{\al, 4m}$
%satisfies a large deviation principle as $\ve \searrow 0$ with a good rate function.
%
%
%
%We assume the coefficient of SDE (\ref{sde_set.eq}) are $C_b^{\infty}$. 
%Note that in the most important theorem,  Theorem \ref{tm.ldp.sgt} (i),
%we will assume the ellipticity condition only at the starting point $a$.

%By Lyons' continuity theorem and 
%the contraction principle (e.g., Thorem 4.2.1, \cite{dzbk}), 
%a large deviation principle for the pushforward of $\mu^{\ve}_{a,a'}$ by 
%the map $\tilde{y}^{\ve}$ immediately follows
%(see Corollary \ref{co.omake} below).
%This can be considered as a special case of Theorem 2.1, \cite{tw}.
%
%
%After normalization, our main theorem, i.e.,
%a large deviation principle for $\{ {\mathbb Q}^{\ve}_{a,a'}\}_{\ve >0}$ 
%under {\bf (H1)}, is easily obtained.

%Let
%$\phi^0 =\phi^0 (h)$ be the unique solution of ODE (\ref{ode_set.eq}). 
%Note that $\phi^0 (h) = a+ \Phi_0 (\tau_{\lambda}(S_2(h)))^1$.
%Under {\bf (H1)}, ${\cal K}^{a,a'}= \{ h \in {\cal H} ~|~  \phi^0 (h)_1 =a'  \}
%$ is a non-empty closed subset of ${\cal H}$.
%In general, however, it can be empty.

%%%%%%%%%%%%%%%%%%%%%%%%
%
%Under {\bf (H1)} , the heat kernel $p^{\ve}_t (a,a')$ exists and positive 
%for all $a, a' \in {\bf R}^n$, $t>0$ and $\ve>0$.
%Hence, the results in the previous section are available.
%

The second assertion 
of Theorem \ref{tm.ldp.sgt} is immediate from the first one,
because the weight of the whole set equals $p^{\ve}_1 (a,a') = {\mathbb E}[ \delta_{a'} (y^{\ve}(1,a)) ]$
and the whole set is both open and closed.
%
%The proof of the first assertion of Theorem \ref{tm.ldp.sgt} will be given in later sections.

As an immediate consequence of 
Theorem \ref{tm.ldp.sgt} {\rm (i)},
we can show the following corollary.
This can be considered as a special case of Theorem 2.1, \cite{tw}.
After normalization, we obtain Theorem \ref{tm.ldp.main}.

%

%Theorem \ref{tm.ldp.sgt}, 
%we can show our main thorem (Theorem \ref{tm.ldp.main}),  i.e., 
%a large deviation principle for the family $\{ {\mathbb Q}^{\ve}_{a,a'}\}_{\ve >0}$
%of pinned diffusion measures.
%The key point is that
%we can use the contraction principle, thanks to Lyons' continuity theorem.
%
%
%

\begin{co}\label{co.omake}
Keep the same notations and assumptions as in Theorem \ref{tm.ldp.sgt}, (i). 
Then, the pushforward measure of $\theta^{\ve}_{a,a'}$ by the map $\tilde{y}^{\ve}$
satisfies a large deviation principle 
with a good rate function $J'(y)$, where
$
J'(y) =\inf\{ \|h\|^2_{{\cal H}}  /2
~|~ h \in {\cal K}^{a,a'}  \mbox{ with }  y = \phi^0 (h) \}
$
if $y = \phi^0 (h)$ for some $h \in {\cal K}^{a,a'}$
and 
$J'(y) = \infty$ if no such $h \in {\cal K}^{a,a'}$ exists.
\end{co}

\vspace{7mm}
\noindent
{\it Proof of Corollary \ref{co.omake} and Theorem \ref{tm.ldp.main}.}\quad
%\Proof
There exists $\al' \in (\al, 1/2)$ and $m \in {\bf N}_+$ such that  $\al' -(4m)^{-1} \ge \al$
and (\ref{eq.amam}) hold.
Then, we can use Theorem \ref{tm.ldp.sgt} above and 
we have a continuous embedding 
$G\Omega^B_{\al', 4m} ({\bf R}^d)  \hookrightarrow  G\Omega^H_{\al} ({\bf R}^d)$
by Proposition \ref{inj.HB.pr}.
Of course, the Young pairing
$\tau_{\lambda}: G\Omega^H_{\al} ({\bf R}^d) \to G\Omega^H_{\al} ({\bf R}^{d+1})$
is continuous.

Let us give a remark on continuity of the Lyons-It\^o map $\Phi_{\ve}$
since it depends on $\ve$, too.
The coefficient $V_0$ for the drift term 
depends on $\ve$ and $C_b^{\infty} ([0,1] \times {\bf R}^n)$.  
Hence, by Taylor expansion, we have 
\[
\sup_{\xi:|\xi| \le R }
| \nabla^i_{(\xi)}  V_0 (\ve, \xi) - \nabla^i_{(\xi)} V_0 (\ve_0, \xi) |
\to 0
\qquad \mbox{ as $\ve \to \ve_0$}
\]
for any $R>0$ and $0 \le i \le 3$.  
Here, $\nabla_{(\xi)} $ stands for the gradient in $\xi$-variable.

Under this kind uniform convergence of the $\ve$-dependent coefficient, 
a slight generalization of Lyons' continuity theorem holds, which claims that
the mapping
\[ 
[0,1] \times G\Omega^H_{\al}   ({\bf R}^{d+1})   
\ni  (\ve, Z)
\mapsto
\Phi_{\ve} (Z)  \in
G\Omega^H_{\al} ({\bf R}^{n})
\]
is continuous.
It immediately follows that the mapping
\[ 
[0,1] \times G\Omega^B_{\al', 4m} ({\bf R}^{d})   
\ni  (\ve, X)
\mapsto
a + \Phi_{\ve} ( \tau_{\lambda} (X) )^1  \in C^{\al -H}_{a} ([0,1], {\bf R}^d)
\]
is continuous, too.
Now we can use Proposition \ref{pr.rep_pin}
and a slight generalization of the contraction principle 
for large deviations (see Lemma 3.9, \cite{ik} for example). 
This completes the proof.

In the same way,
the proof of Theorem \ref{tm.ldp.main} follows from Proposition \ref{pr.rep_pin} 
and Theorem \ref{tm.ldp.sgt}, (ii).
\toy

Before we prove Theorem \ref{tm.ldp.sgt} {\rm (i)},
let us recall that a Schilder-type large deviation for the laws of scaled Brownian rough path $\ve W$ holds 
on $G\Omega^{H}_{\al'} ({\bf R}^d)$ as $\ve \searrow 0$
for any $1/3< \al' <1/2$ with a good rate function $I_{Sch}$
(see Friz-Victoir \cite{fv1}).
Here, $I_{Sch} (X) = \|h\|^2_{{\cal H}}/2 $ if $X=S_2(h)$ for some $h \in {\cal H}$
and $I_{Sch} (X) = \infty$ if otherwise.
For $\al <\al' <1/2$, there is a continuous embedding 
$G\Omega^{H}_{\al'}  ({\bf R}^d) \hookrightarrow G\Omega^{B}_{\al, 4m} ({\bf R}^d)$.
Hence, by the contraction principle,
the Schilder-type large deviation also holds on $G\Omega^{B}_{\al, 4m} ({\bf R}^d)$
and 
$I_{Sch} $ is still good on $G\Omega^{B}_{\al, 4m} ({\bf R}^d)$.
We can easily verify from this that
$I$ defined in Section 2
 is lower semicontinuous and good (i.e. the level set $\{ X| I(X) \le c\}$ is compact 
for all $0 \le c < \infty$).  So is $\hat{I}$.
%
%

%\begin{re}\label{re.anoI}
%Let $F$ be any continuous map from $G\Omega^H_{\al} ({\bf R}^{d})$ 
%to a Hausdorff topological space $(1/3 <\al < 1/2)$.
%Then, from the proof of the above corollary, 
%it is obvious that the family  $\{ F_* \mu^{\ve}_{a,a'} \}_{\ve >0}$ 
%of image measures satisfies a large deviation principle
%as $\ve \searrow 0$.
%In particular, we may take another Lyons-It\^o map  as $F$.  
%This is worth mentioning since Theorem 2.1 \cite{tw} is formulated in that way.
%\end{re}

%%%%%%%%%%%%%%%%%%%%%%%%%%%%%%%%%%%%%%%%%%%%%%%%%%%%%%%%%%%%%%%%%%%%%%%%%%%%%%%%%%%%%%%%%%%%%%
%%%%%%%%%%%%%%%%%%%%%%%%%%%%%%%%%%%%%%%%%%%%%%%%%%%%%%%%%%%%%%%%%%%%%%%%%%%%%%%%%%%%%%
%%%%%%%%%%%%%%%%%%%%%%%%%%%%%%%%%%%%%%%%%%%%%%%%%%%%%%%%%%%%%%%%%%%%%%%%%%%%%%%%%%%%
%%  Section    Proof :  Lower Bound
%%%%%%%%%%%%%%%%%%%%%%%%%%%%%%%%%%%%%%%%%%%%%%%%%%%%%%%%%%%%%%%%%%%%%%%%%%%%%%%%%%%%%%
%%%%%%%%%%%%%%%%%%%%%%%%%%%%%%%%%%%%%%%%%%%%%%%%%%%%%%%%%%%%%%%%%%%%%%%%%%%%%%%%%%%%%%
%%%%%%%%%%%%%%%%%%%%%%%%%%%%%%%%%%%%%%%%%%%%%%%%%%%%%%%%%%%%%%%%%%%%%%%%%%%

\subsection{Proof of Theorem \ref{tm.ldp.sgt} {\rm (i)}: Lower estimate}
\label{se.low}
The aim of this subsection is to prove the lower estimate
in Theorem \ref{tm.ldp.sgt} {\rm (i)}.

Before we start our proof, let us rewrite SDE (\ref{sde_set.eq})
in matrix form.
Set $\sigma : {\bf R}^n \to {\rm Mat} (n, d)$ by $[V_1, \ldots, V_d]$,
that is, the $j$th column vector of $\sigma$ is $V_j$.
We also write $V_0 (\ve, \xi) =b(\ve, \xi)$.
Then, SDE (\ref{sde_set.eq}) can be rewritten as follows;
\begin{equation}\label{sde.def}
dy^{\ve}_t  =  \sigma( y^{\ve}_t) \circ \ve dw_t + b(\ve,  y^{\ve}_t)dt,
\qquad\qquad
 y^{\ve}_0 = a.
\end{equation}
This matrix form (\ref{sde.def}) is sometimes more convenient than  (\ref{sde_set.eq}).
In the same way, ODE (\ref{ode_set.eq}) is equivalent to the following ODE;
\begin{equation}
%\label{sde.def}
d \phi^{0}_t  = \sigma( \phi^{0}_t)  dh_t + b(0,  \phi^{0}_t)dt,
\qquad\qquad
 \phi^{0}_0 = a.
\nn
\end{equation}

%%%%%%%%%%%%%%%%%%%%%%%%%%%%%%%%%%%%%%%%%%%%%%%%%%%%%%%%%%%%%%%%%%%%%%%%%
%%%%%%%%%%%%%%%%%%%%%%%%%%%%%%%%%%%%%%%%%%%%%%%%%%%%%%%%%%%%%%%%%%%%%%%%%

%%%%%%%%%%%%%%%%%%%%%%%%%%%%%   Shifted  SDE  %%%%%%%%%%%%%%%%%%%%%%%%%%

Let $(y^{\ve}_t (w))_{ 0 \le t \le 1}$ be the solution of SDE (\ref{sde.def}).
For $h \in {\cal H}$, we consider the Cameron-Martin shift of $y^{\ve}$.
%
%
%we set $y^{\ve, h}_t (w) = y^{\ve}_t (w+ h/\ve)$.
%Then, it satisfies the following SDE;
%
%
\begin{equation}\label{shift_sde.def}
dy^{\ve, h}_t  = \sigma( y^{\ve, h}_t) [\circ \ve  dw_t + dh_t]+ b(\ve,  y^{\ve, h}_t)dt,
\qquad\qquad
 y^{\ve,h}_0 = a.
\end{equation}
%
%
%
%An equivalent expression is as follows;
%\begin{equation}\label{shift_sde2.def}
%dy^{\ve, h}_t  =
%  \sum_{i=1}^d  V_i( y^{\ve,h}_t) [ \circ \ve  dw_t^i + dh_t^i] +  V_0 (\ve,  y^{\ve,h}_t)dt, 
%\qquad\qquad
% y^{\ve, h}_0 = a.
%\end{equation}
%
%
An $\infty$-quasi-continuous version of $y^{\ve, h}$ is given by 
$a + \Phi_{\ve} ( \tau_{\lambda} ( S_2 (\ve w +h)) )$ 
by Lemma \ref{lm.Zjimei} and Proposition \ref{pr.MalNua}.

%%%%%%%%%%%%%%

Non-degeneracy of the Malliavin covariance matrix of the solution  of SDE 
under ellipticity condition is in Section 6.1.5, Shigekawa \cite{sh} for example. 
By repeating the same proof with carefully taking $\ve$-dependency into account, 
we can verify that $(y^{\ve, h}_1 - a' ) /\ve$ 
is uniformly non-degenerate in the sense of Malliavin as $\ve \searrow 0$.
(This fact was already shown in \cite{wa} when the drift term is 
of the form $b(\ve, \xi)= \ve^2 \hat{b}(\xi)$ for an $\ve$-independent 
vector field $\hat{b}$.
The definition of uniform non-degeneracy is in p. 10, \cite{wa}.
When we try to extend our method to the hypoelliptic case, 
this part becomes very difficult.)

Moreover, in the same way as in Watanabe \cite{wa}, we have the following
asymptotics for $y^{\ve, h}_1$ as $\ve \searrow 0$;
\[
y^{\ve, h}_1 = \phi^0 (h)_1 + \ve  \phi^1 ( \,\cdot\, ;h)_1  + O(\ve^2)
\qquad
\mbox{ in ${\bf D}_{\infty} ({\bf R}^n)$.}
\]
Watanabe obtained the expansion up to any order 
when the drift term is of the form $b(\ve, \xi)= \ve^2 \hat{b}(\xi)$. 
But we only need an expansion 
up to second order. Modification for the case of general $b(\ve, \xi)$ is easy.
Here, $\phi^1_t = \phi^1 (w; h)_t$ is the unique solution of the following equation;
\[
d \phi^1_t - \nabla \sigma ( \phi^0_t) \la \phi^1_t, dh_t\ra 
- \nabla_{(\xi)} b (0,  \phi^0_t)  \la \phi^1_t \ra dt
=
\sigma ( \phi^0_t ) dw_t + \partial_{\ve} b (0,  \phi^0_t) dt,
\quad
\phi^1_0 =0.
\]
Here, $ \nabla_{(\xi)}$ denotes the (partial) gradient with respect to the second variable.
$\partial_{\ve}$ should be understood in a similar way.

Now we give an explicit expression for $\phi^1_t$.
Let $M_t= M(h)_t$ be a unique solution of the following $n \times n$ matrix-valued ODE;
\[
dM_t =  [ \nabla \sigma ( \phi^0_t) \la \bullet , dh_t\ra 
+ \nabla_{(\xi)} b (0,  \phi^0_t)  \la \bullet \ra dt] \cdot M_t,
\qquad
M_0 ={\rm Id}_{n}.
\]
Then, $M_t$ is invertible for all $t$ and we have 
\[
 \phi^1_t  = M_t  \int_0^t  M_s^{-1} [ \sigma ( \phi^0_s ) dw_s + \partial_{\ve} b (0,  \phi^0_s) ds].
\]
Hence, as a functional in $w$,  $\phi^1_1$ is an ${\bf R}^n$-valued Gaussian random variable
with mean $M_1  \int_0^1  M_s^{-1}  \nabla_{\ve} b (0,  \phi^0_s) ds$ and 
 covariance 
$M_1  \int_0^1  M_s^{-1}  \sigma ( \phi^0_s ) \sigma ( \phi^0_s )^T  (M_s^{-1})^T ds  M_1^T$.
Since we assumed ellipticity at the initial point $a$,
$\sigma \sigma^T$ is positive and bounded away from $0$ near $a$.
Thereore, this covariance matrix is non-degenerate for any $h$.
In particular,  ${\mathbb E} [\delta_0 (\phi^1_1) ]>0$.

Let $h \in {\cal K}^{a,a'}$.
By (a special case of) Theorem 2.3, \cite{wa} or Theorem 9.4, \cite{iwbk} 
and the uniform non-degeneracy, 
we have
$$
\lim_{ \ve \searrow 0 }
\delta_0 \bigl( \frac{y^{\ve, h}_1 - a' }{\ve} \bigr) = \delta_0( \phi^1_1)
\qquad
\mbox{ in ${\bf D}_{q, -r} ({\bf R}^n)$}
$$
for some $q \in (1, \infty)$ and some $r \in {\bf N}$.

%%%%

\begin{lm}\label{lm.dnsK}
${\cal K}^{a,a'} \cap {\cal W}^*$ is dense in ${\cal K}^{a,a'}$.
In other words,
for any $h \in {\cal K}^{a,a'}$,  there exist $h_j\in {\cal K}^{a,a'} ~(j \in {\bf N}_+)$
such that $\la h_j, \,\cdot\, \ra_{{\cal H}} \in  {\cal W}^*$ for all $j$ and 
 $\lim_{j \to \infty} \| h_j -h  \|_{{\cal H}} =0$.
\end{lm}

\Proof
Set $F(h) =\phi^0(h)_1$.
Then, $F: {\cal H} \to {\bf R}^n$ is Fr\'echet-$C^1$ and its Fr\'echet derivative $D F$
is given by 
\[
D_k F(h)
=
\la D F(h), k \ra_{{\cal H}}
=
M_1  \int_0^1  M_s^{-1}  \sigma ( \phi^0(h)_s )  k'_s ds. 
\]
Surjectivity of the linear map $DF (h): {\cal H} \to {\bf R}^n$
is equivalent to non-degeneracy of the "deterministic Malliavin covariance" matrix
$$
(\la D F^i(h),   D F^j(h)\ra_{{\cal H}} )_{1 \le i,j \le n}
=
M_1  \int_0^1  M_s^{-1}  \sigma ( \phi^0(h)_s )\sigma ( \phi^0(h)_s )^T (M_s^{-1})^T ds  M_1^T,
$$
which has already been shown.
Hence, we can use Lemma \ref{lm.inver} below with ${\cal K}={\cal H}$ and ${\cal L}={\cal W}^*$.
\QED

%%%%%%%%%%%%%%%%%%%%%%%%%%%%%%%

\begin{lm}\label{lm.inver}
Let ${\cal K}$ be a real Hilbert space and $\xi \in {\cal K}$.
Assume that {\rm (i)}~ $F$ is an ${\bf R}^n$-valued Fr\'echet-$C^1$ map 
 defined on a neighborhood of $\xi$ with a bounded derivative $D F$
and {\rm (ii)}~$DF (\xi): {\cal K} \to {\bf R}^n$ is a surjective linear map.
Let ${\cal L}$ be a real Banach space which is continuously and densely 
embedded in ${\cal K}$.
Then, there exists $\xi_j \in {\cal L}~(j=1,2,\ldots)$ such that 
$\lim_{j\to \infty} \| \xi_j- \xi \|_{\cal K} =0$ and $F(\xi_j) = F(\xi)$ for all $j$.
\end{lm}

\Proof
Without loss of generality, we may assume $F(\xi)=0$.
It is sufficient to show that, for sufficiently small $r>0$,
there exists
$\eta \in {\cal L} \cap B_{{\cal K}}(\xi, r)$
such that $F(\eta)$=0, 
where $B_{{\cal K}}(\xi, r)$ is the open ${\cal K}$-ball of radius $r$, centered at $\xi$.
There exists an $n$-dimensional subspace $V \subset {\cal K}$
such that $\nabla F(\xi)|_V : V \to {\bf R}^n$ is a linear bijection.
Hence, we may apply inverse function theorem at $\xi$ to $F|_{\xi +V}$,
which is the restriction of $F$ onto the affine subspace $\xi +V$.

First, we consider the case for $n=1$.
For any sufficiently small $r>0$,
  there exist $\chi_0, \chi_1 \in (\xi +V) \cap B_{{\cal K}}(\xi, r)$
such that $F(\chi_0) >0$ and $F(\chi_1) <0$.
Since ${\cal L}$ is dense in ${\cal K}$, we can find 
$\eta_0, \eta_1 \in {\cal L} \cap B_{{\cal K}}(\xi, r)$
such that $F(\eta_0) >0$ and $F(\eta_1) <0$.
Then, the line segment $\overline{\eta_0 \eta_1}
 =\{\tau \eta_0 + (1-\tau)\eta_1 ~|~ 0 \le \tau \le 1\}$ 
 is still in ${\cal L} \cap B_{{\cal K}}(\xi, r)$.
Then, by the intermediate value theorem, there exists $\tau' \in (0,1)$
such that $F( \tau' \eta_0 + (1-\tau')\eta_1)=0$.
Thus, we have shown the case $n=1$.

Next, we consider the case for $n=2$.
For $0<r<r'$, set $B_{{\bf R}^n}(0, r)= \{ a \in {\bf R}^n ~|~ |a| < r \}$,
$S^1_r =\{ a \in {\bf R}^n ~|~ |a| = r \}$,
and
$A_{r,r'} =\{ a \in {\bf R}^n ~|~ r < |a| < r' \}$.
For sufficiently small $r>0$, $F ((\xi +V) \cap B_{{\cal K}}(\xi, r))$
is an open neighborhood of $0$ and 
$f:=F|_{ (\xi +V) \cap B_{{\cal K}}(\xi, r)}$ is a diffeomorphism.
There exists $\rho >0$ such that 
$B_{{\bf R}^n}(0, 3\rho) \subset F ((\xi +V) \cap B_{{\cal K}}(\xi, r))$.
Then, there exist $N \in {\bf N}_+$ and $\xi_0, \ldots, \xi_N \in S^1_{\rho}$ such that
the following conditions are satisfied;

\begin{enumerate}
\item
 $\xi_i$ and $\xi_{i+1}$
are adjacent in counter-clockwise order (with $\xi_{N+1} =\xi_0$),
\item
the arc 
%$\overarc{\xi_i \xi_{i+1}}$
%$\stackrel{\frown}{\xi_i \xi_{i+1}}$ 
between $\xi_i$ and $\xi_{i+1}$ 
can be continuously deformed inside $A_{\rho/2 ,2\rho}$,
with the end points being fixed, so that the image of the arc by $f^{-1}$ is deformed 
to the line segment $\overline{ f^{-1}(\xi_i) f^{-1}(\xi_{i+1})}$.
\end{enumerate}

Then, there exists $\eta_0, \ldots, \eta_N \in {\cal L} \cap B_{{\cal K}}(\xi, r)$
such that the following conditions are satisfied;

\begin{enumerate}
\item
there exists a continuous deformation of 
$\cup_{i=0}^N \overline{ f^{-1}(\xi_i) f^{-1}(\xi_{i+1})}$ inside $B_{{\cal K}}(\xi, r)$,
 which deforms
each  $\overline{ f^{-1}(\xi_i) f^{-1}(\xi_{i+1})}$
to  $\overline{\eta_i \eta_{i+1}}$.
\item
the image of the above continuous deformation by $F$ stays inside $A_{\rho/3 ,3\rho}$.
\end{enumerate}

Obviously, $\cup_{i=0}^N \overline{\eta_i \eta_{i+1}}$ is a subset of 
the following (possibly degenerate) simplex
\[
\triangle :=\{ \eta_0 + \sum_{i=1}^N \tau_i (\eta_i -\eta_0) 
~|~ \mbox{$\tau_i \ge 0$  for all $i$ and $\sum_{i=1}^N \tau_i \le 1$} \}.
\]
This is the convex hull of $\eta_0, \ldots, \eta_N$.
Hence,
$\triangle \subset B_{{\cal K}}(\xi, r) \cap {\cal L}$.
It is easy to see that $\triangle$ admits a continuous deformation to a single point set $\{\eta_0\}$
without getting out of $\triangle$.

Combining them all, we have shown that $S^1_{\rho}$ is continuously 
deformed in ${\bf R}^2$ to a single point set $\{F(\eta_0)\}$.
More precisely, there exists a continuous map $\phi :[0,1] \times S^1_{\rho} \to {\bf R}^2$
such that $\phi(0, \,\cdot\,) ={\rm Id}_{ S^{1}_{\rho}}$ and $\phi (1, \,\cdot\,)$ is constant.
By Lemma \ref{lm.mrys} below, there exists 
$(s, a) \in [0,1] \times S^{1}_{\rho}$ such that $\phi(s,a)=0$.
By way of construction of the deformation map, 
there must be $\eta \in \triangle$ such that $F(\eta)=\phi(s,a)=0$.
This proves the case $n=2$.
The higher dimensional cases can be shown in the same way.
\QED

The following topological lemma can be regarded as a higher dimensional version of 
the intermediate value theorem. 

\begin{lm}\label{lm.mrys}
Let $B^n = \{ a \in {\bf R}^n ~|~ |a| \le 1 \}$ be the $n$-dimensional closed ball
and let
$S^{n-1} = \{ a \in {\bf R}^n ~|~ |a| =1 \}$ be the $(n-1)$-dimensional sphere
($n \ge 2$). 
\\
\noindent
{\bf (i)}~ Let $f : [0,1] \times S^{n- 1} \to {\bf R}^n$ be a continuous map such that 
$f(0, \,\cdot\,) = {\rm Id}_{ S^{n- 1}}$ and $f(1, \,\cdot\,)$ is a constant map.
Then, there exists $(s, a) \in [0,1] \times S^{k- 1}$ such that $f(s,a)=0$.
\\
\noindent
{\bf (ii)}~
Let $g : B^{n} \to {\bf R}^n$ be a continuous map such that 
$g|_{ S^{n- 1}} = {\rm Id}_{ S^{n- 1}}$.
Then, there exists $a\in B^{n}$ such that $g(a)=0$.
\end{lm}

\Proof
Since $f(1, \,\cdot\,)$ is constant, there exists a continuous map
$g: B^{n} \to {\bf R}^n$ such that 
$g|_{ S^{n- 1}} = {\rm Id}_{ S^{n- 1}}$ and ${\rm Im} f = {\rm Im} g$.
Hence, it is sufficient to prove the second assertion.
Assume that $g^{-1} (0) =  \emptyset$.
Then $g$ is continuous from $B^n$ to ${\bf R}^n \setminus \{0\}$
and we have the following commutative diagram:

$$
\begin{CD}
\{0\}=H_{n-1} (B^n) @> g_*  >>   H_{n-1} ({\bf R}^n \setminus \{0\}) ={\bf Z}\\
@A\iota_* AA    @AA \iota_* A \\
{\bf Z}=H_{n-1} (S^{n- 1})  @>  ( g|_{ S^{n- 1}})_*  >>  H_{n-1} (S^{n- 1})={\bf Z}
\end{CD}
$$
where $\iota$ stands for the injections.
$g_* \circ (\iota: S^{n-1}  \hookrightarrow B^n)_*$ is the zero map.
On the other hand, 
$(\iota:S^{n-1} \hookrightarrow {\bf R}^n \setminus \{0\})_* \circ ( g|_{ S^{n- 1}})_*
=(\iota:S^{n-1} \hookrightarrow {\bf R}^n \setminus \{0\})_* \circ {\rm id}_{ H_{n-1} (S^{n- 1})}$
is an isomorphism, due to homotopy equivalence.
This is a contradiction.
Therefore, $g^{-1} (0) \neq \emptyset$.
\QED

%%%%%%%%%%%%%%%%%%%%%%%%%%%%%%

The lower estimate in Theorem \ref{tm.ldp.sgt} {\rm (i)} is almost immediate
from Lemma \ref{lm.dnsK} above and  the following proposition.
\begin{pr}\label{pr.lowball}
Let $A \subset G\Omega^B_{\al, 4m} ({\bf R}^{d})$ be open and 
suppose that $S_2(h) \in A$
and $\la h, \,\cdot\, \ra_{{\cal H}} \in  {\cal W}^* \cap {\cal K}^{a,a'}$.  
Then, we have $- \|  h \|^2_{{\cal H}} /2
\le
 \liminf_{\ve \searrow 0 } \ve^2 \log \mu^{\ve}_{a,a'} (A)$.
\end{pr}

\Proof
For $R>0$, we set 
$
B_{R}=\{ X \in G\Omega^B_{  \al, 4m} ({\bf R}^{d}) ~|~ \|X^i \|_{i \al, 4m/i -B}  <R^{i}
\quad (i=1,2) \}
$
and 
$\hat{B}_{R} (H)= T_h ( B_{R} )$, 
where $T_h$ is the Young translation by $h$ on $G\Omega^B_{\al, 4m} ({\bf R}^{d})$.
Then, $\{\hat{B}_{R} (H) ~|~ R>0\}$ forms a fundamental system of open neighborhood around $H =S_2(h)$.
Since $A$ is open, there exists $R_0 >0$ such that 
$
\hat{B}_{R_0} (H) \subset A.
$
We will estimate the weight of $\hat{B}_{R_0}  (H)$ instead of that of $A$.

By Cameron-Martin formula, it holds that,
for any $F \in {\bf D}_{\infty}$, 
\begin{align}
{\mathbb E} [ F  \delta_{a'} (y^{\ve}_1) ]
&=
{\mathbb E} [ \exp (- \frac{ \la h, w\ra }{\ve} - \frac{ \| h\|^2_{{\cal H}} }{2\ve^2} ) 
 F( w +\frac{h}{\ve})  \delta_{a'} (y^{\ve, h}_1) ]
\nn\\
&=
e^{ - \| h\|^2_{{\cal H}}/2\ve^2 }  \ve^{-n} 
{\mathbb E} [ e^{ - \la h, w\ra /\ve } 
 F( w+\frac{h}{\ve})  \delta_{0} ( \frac{ y^{\ve, h}_1  -a'}{\ve}) ].
\nn
\end{align}
Here, we used the fact that $\delta_0 (\ve \xi) = \ve^{-n} \delta_0 (\xi)$.

We denote by $\nu^{\ve} = \nu^{\ve}_{a,a'}$ the Borel measure corresponding to 
$\delta_{0} ( (y^{\ve, h}_1  -a' )/\ve )$ via Sugita's theorem (Proposition \ref{pr.sugita}).
Then,  we have
\begin{align}
\mu^{\ve}_{a,a'} ( \hat{B}_R (H) ) 
&=
\int_{{\cal W}}     I_{\hat{B}_R (H) } (\ve W)     \,  \theta^{\ve}_{a,a'}  (dw)
\nn\\
&=
e^{ - \| h\|^2_{{\cal H}}/ 2\ve^2}  \ve^{-n} 
\int_{{\cal W}}  e^{- \la h, w\ra /\ve } 
   I_{ B_{R/ \ve} }  (W) \nu^{\ve} (dw).
\nn
\end{align}
for any $0 <R<R_0$.
Noting that $|\la h, w\ra /\ve| \le \|h\|_{{\cal W}^*} \| w/\ve \|_{{\cal W}} 
\le C\|h\|_{{\cal W}^*} R/ \ve^2$,
where $C := \sup_{w \neq 0} \| w \|_{\infty} / \| w \|_{\al, 4m -B}$, 
we see that 
\[
\mu^{\ve}_{a,a'} ( \hat{B}_R (H) ) 
\ge
e^{ - \| h\|^2_{{\cal H}}/(2\ve^2) }  \ve^{-n}   e^{- C\|h\|_{{\cal W}^*} R /\ve^2 }
 \nu^{\ve} \bigl(   \{ \|W^i  \|_{i \al, 4m/i -B}^{1/i}  <R /\ve
\,\, (i=1,2)\}  \bigr).
\]

Hence, it is sufficient to show that, for each fixed $R$, 
\begin{equation}\label{sufflow.eq}
\lim_{\ve \searrow 0} \nu^{\ve} \bigl(   \{ \|W^i  \|_{i \al, 4m/i -B}^{1/i}  <R /\ve
\,\, (i=1,2)\}  \bigr)
=
\nu^{\infty} ( {\cal W}) =
 {\mathbb E}[  \delta_0 ( \phi^1_1 )] >0,
\end{equation}
where $\nu^{\infty}$ stands for the Borel measure corresponding to 
$\delta_{0} (  \phi^1_1 )$ via Sugita's theorem.
Indeed,  (\ref{sufflow.eq}) implies that
\begin{align}
\liminf_{  \ve \searrow 0}  \ve^2 \log \mu^{\ve}_{a,a'} ( \hat{B}_{R_0} (H) ) 
\ge
\liminf_{  \ve \searrow 0}  \ve^2 \log \mu^{\ve}_{a,a'} ( \hat{B}_R (H) ) 
%\nn\\&
\ge
-  \frac{  \| h\|^2_{{\cal H}}   }{2 }     -   C\|h\|_{{\cal W}^*} R 
\nn
\end{align}
Letting $R \searrow 0$, we have
$
\liminf_{  \ve \searrow 0}  \ve^2 \log \mu^{\ve}_{a,a'} ( \hat{B}_{R_0} (H) ) 
\ge
- \| h\|^2_{{\cal H}}/2. 
$

It remains to prove (\ref{sufflow.eq}).
Since $\delta_{0} ( (y^{\ve, h}_1  -a' )/\ve )   \to \delta_0 ( \phi^1_1 )$ 
in ${\bf D}_{q', -r}$ for some $q' \in (1, \infty)$ and $r \ge 0$, 
we have 
$\nu^{\ve} ( {\cal W}) -\nu^{\infty} ( {\cal W}) 
=  
{\mathbb E}[  1 \cdot  \{  \delta_{0} ( (y^{\ve, h}_1  -a' )/\ve ) -  \delta_0 ( \phi^1_1 ) \} ]  \to 0$
as $\ve \searrow 0$.
Now, it is sufficient to show that
\begin{align} 
\lefteqn{
\nu^{\ve} \bigl(   \{ \|W^i  \|_{i \al, 4m/i -B}^{1/i}  <R /\ve
\,\, (i=1,2)\}^c  \bigr)
\le
\sum_{i=1,2} 
\nu^{\ve} \bigl(   \{ \|W^i  \|_{i \al, 4m/i -B}^{1/i}  \ge R /\ve\} \bigr)   
}
\nn\\
&\le 
\sum_{i=1,2}  \| \delta_{0} ( (y^{\ve, h}_1  -a' )/\ve ) \|_{q', -r} 
\cdot
{\rm c}_{q,r} \bigl(   \{ \|W^i  \|_{i \al, 4m/i -B}^{1/i}  \ge R /\ve\} \bigr)
\to 0
\nn
\end{align}
as $\ve \searrow 0$.
Here, $1/q +1/q' =1$.
Since $\| \delta_{0} ( (y^{\ve, h}_1  -a' )/\ve ) \|_{q', -r}$ is bounded in $\ve$, 
the problem is reduced to showing 
\[
\lim_{\ve \searrow 0}
{\rm c}_{q,r} \bigl(   \{ \|W^i  \|_{i \al, 4m/i -B}^{1/i}  \ge R /\ve\} \bigr)   =0  
\qquad
\qquad
(i=1,2).
\]
This will be shown in the next lemma below.
\QED

%%%%%%%%%%%%%%%%%%%%%%%%%%%%%%%%%%%%%%%%%%%%%%%%%%%%%%%%%%%%%%%%%%%%%%%%%%%%%%%%%%%%%
%%%%%%%%%%%%%%%%%%%%%%%%%%%%%%  Lemma    cap Gauss decay  %%%%%%%%%%%%%%%%%%%%%%%%%%%%%
%%%%%%%%%%%%%%%%%%%%%%%%%%%%%%%%%%%%%%%%%%%%%%%%%%%%%%%%%%%%%%%%%%%%%%%%%%%%%%%%%%%%%

On an abstract Wiener space, weight of the complement set of the large ball 
admits Gaussian decay.
This is called a large deviation estimate.
A similar fact holds in rough path settings, too (for instance, see \cite{fo}).
In this paper, however, we need a large deviation estimate for capacities, not for measures.

\begin{lm}\label{lm.capdcy}
Let $\al$ and $m$ be as in Theorem \ref{tm.ldp.sgt}.
For any $q \in (1, \infty)$ and $r \in {\bf N}$, there exist $c>0$ and $R_1>0$ such that
\[
{\rm c}_{q,r} \bigl(  \{ w \in {\cal W} ~|~   \|W^i  \|_{i \al, 4m/i -B}^{1/i}  \ge R \} \bigr) 
\le 
e^{ - c R^2}
\qquad
\mbox{for any  $R \ge R_1$ and $i=1,2$}.  
\]
The constants $c$ and $R_1$ may depend on $\al, m, q, r$, but not on $R$.
\end{lm}

\Proof
First, let us check that there exists $c_1 >0$ such that 
\begin{equation}\label{isop.ineq}
{\mathbb P} \bigl(  \{ w \in {\cal W} ~|~   \|W^i  \|_{i \al, 4m/i -B}^{1/i}  \ge R \} \bigr) 
\le 
e^{ - c_1 R^2}
\end{equation}
holds for sufficiently large $R>0$.
On the H\"older geometric rough path space $G\Omega^H_{\al} ({\bf R}^d)$,
with any $1/3 <\al < 1/2$,
this type of Gaussian decay is well-known
(See Friz and Oberhauser \cite{fo}).
Next, take $\al' \in (\al, 1/2)$. Then, due to the continuous embedding 
$G\Omega^H_{\al'} ({\bf R}^d)  \hookrightarrow G\Omega^B_{\al, 4m} ({\bf R}^d)$,
we can easily see (\ref{isop.ineq}) holds on the Besov geometric rough path space, too.

Take $\chi \in C_b^{\infty} ({\bf R} \to {\bf R})$ such that 
$\chi =0$ on $(-\infty, 0]$,  $\chi =1$ on $[1,\infty)$, and $\chi$ is non-decreasing.
Set $G(w) = \|W^1 \|_{ \al, 4m -B}^{4m}$.
Then, $G \in {\bf D}_{\infty}$ and $\chi (G -R^{4m} +1) 
\in {\cal F}_{q,r} ( \{\|W^1 \|_{\al, 4m -B}  \ge R  \} )$
for any $q, r$ in the sense of (\ref{eq.libr.def}), since $w \mapsto W$ is 
$\infty$-quasi-continuous by Corollary \ref{co.qcontW}.
By (\ref{lib.pot.ineq}), we have
$$
{\rm c}_{q,r} \bigl(  \{\|W^1 \|_{\al, 4m -B}  \ge R  \} \bigr) \le \| \chi (G -R^{4m} +1)  \|_{q,r}.
$$
So, we have only to estimate the Sobolev norm on the right hand side.

First we calculate the $L^q$-norm.
\[
\| \chi (G -R^{4m} +1)  \|_{q} 
\le {\mathbb P} \bigl(  \{    \|W^1 \|_{\al, 4m-B}  \ge (R^{4m} -1)^{1/4m} \} \bigr)^{1/q}
\le 
e^{ - (c_1 /2q)R^2}
\]
for sufficiently large $R$.
It is easy to see that $D [\chi (G -R^{4m} +1) ]= \chi^{\prime} (G -R^{4m} +1) DG$.
Hence, 
\begin{align}
\| D [\chi (G -R^{4m} +1) ]\|_{q}  
&\le 
\| \chi^{\prime} \|_{\infty}    \bigl\| I_{  \{ G -R^{4m} +1 \ge 0\} }  \|DG\|_{{\cal H}}  \bigr\|_{q} 
\nn\\
&\le
\| \chi^{\prime} \|_{\infty}  \|DG\|_{2q}
 {\mathbb P} \bigl(  \{    \|W^1 \|_{\al, 4m-B}  \ge (R^{4m} -1)^{1/4m} \} \bigr)^{1/2q}
\nn\\
&\le
\| \chi^{\prime} \|_{\infty}  \|DG  \|_{2q}   \cdot
e^{ - (c_1 /4q)R^2}
\nn
\end{align}
for sufficiently large $R$.
From this and Meyer's equivalence,  there exists $c_2 =c_2(\al, m,q)>0$ such that 
\[
\| \chi (G -R^{4m} +1)  \|_{1, q}  \le e^{ - c_2R^2}
\]
for sufficiently large $R$.
In the same way, we can estimate $D^k [\chi (G-R^{4m} +1)]$ for any $k \in {\bf R}_+$
and prove the lemma for $i=1$.

For $i=2$,  just consider 
$$\tilde{G}(w) = \|W^2 \|_{ 2\al, 2m -B}^{2m}
\quad \mbox{ and }\quad
\chi ( \tilde{G} -R^{4m} +1) \in {\cal F}_{q,r} ( \{\|W^2 \|_{2\al, 2m -B}^{1/2}  \ge R  \} .
$$
Then, the same argument works for this case, too. 
\QED

%%%%%%%%%%%%%%%%%%%%%%%%%%%%%%%%%%%%%%%%%%%%%%%%%%%%%%%%%%%%%%%%%%%%%%%%%%%%%%%%%%%%%%%%%%%%%%
%%%%%%%%%%%%%%%%%%%%%%%%%%%%%%%%%%%%%%%%%%%%%%%%%%%%%%%%%%%%%%%%%%%%%%%%%%%%%%%%%%%%%%
%%%%%%%%%%%%%%%%%%%%%%%%%%%%%%%%%%%%%%%%%%%%%%%%%%%%%%%%%%%%%%%%%%%%%%%%%%%%%%%%%%%%
%%  Section    Proof :  upper Bound
%%%%%%%%%%%%%%%%%%%%%%%%%%%%%%%%%%%%%%%%%%%%%%%%%%%%%%%%%%%%%%%%%%%%%%%%%%%%%%%%%%%%%%
%%%%%%%%%%%%%%%%%%%%%%%%%%%%%%%%%%%%%%%%%%%%%%%%%%%%%%%%%%%%%%%%%%%%%%%%%%%%%%%%%%%%%%
%%%%%%%%%%%%%%%%%%%%%%%%%%%%%%%%%%%%%%%%%%%%%%%%%%%%%%%%%%%%%%%%%%%%%%%%%%%

\subsection{Proof of Theorem \ref{tm.ldp.sgt} {\rm (i)}: Uppeer estimate}
\label{se.up}
The aim of this subsection is to prove the upper estimate
in Theorem \ref{tm.ldp.sgt} {\rm (i)}.
In this section, we will often use the following fact;
For $f, g: (0,1] \to [0, \infty)$, it holds that 
$
\limsup_{\ve \searrow 0} \ve^2 \log ( f_{\ve} + g_{\ve}  )  
\le
[ \limsup_{\ve \searrow 0} \ve^2 \log  f_{\ve}  ] \vee 
 [ \limsup_{\ve \searrow 0} \ve^2 \log  g_{\ve} ].
$

%%%%%%%%%%%%%%%%%%%%%

{\bf [Step 1]}~
We divide the proof into three steps. 
The first step is to show that 
\begin{equation}\label{limsmball.ineq}
\lim_{R \searrow 0}  \limsup_{\ve \searrow 0}  \ve^2 \log \mu_{a,a'}^{\ve} ( B_{R} (X) )  \le   - I(X),
\qquad 
X \in G\Omega^B_{\al, 4m}  ({\bf R}^d),
\end{equation}
where 
\[
B_{R} (X)=\{ Y \in G\Omega^B_{  \al, 4m} ({\bf R}^{d}) ~|~ \|Y^i -X^i \|_{i \al, 4m/i -B}  <R^{i}
\quad (i=1,2) \}.
\]

First, we consider the case $a+ \Phi_{0} ( \tau_{\lambda} (X))^1_{0,1} \neq a'$.
For simplicity we write $\tilde{a} = a+ \Phi_{0} ( \tau_{\lambda} (X))^1_{0,1}$.
By continuity of Lyons-It\^o map, 
there exist $\ve_0 >0$ and $R >0$ such that 
$|a+\Phi_{\ve} ( \tau_{\lambda} (Z))^1_{0,1}   -\tilde{a} | \le |a' - \tilde{a}|/3$
holds for all $0 \le \ve \le \ve_0$ and $Z \in B_{R} (X)$.
If we have $\mu_{a,a'}^{\ve} ( B_{R} (X) ) =0$ for such $\ve$ and $R$,
then (\ref{limsmball.ineq}) immediately follows for this case.

Let us verify the fact $\mu_{a,a'}^{\ve} ( B_{R} (X) ) =0$ as follows.
Let $\psi: {\bf R} \to [0,1]$ be a smooth even function such that
$\psi = 1$ on $[0,1]$ and $\psi  = 0$ on $[2, \infty)$ and  non-increasing on $[1,2]$.
Set $\psi (| ( y^{\ve}_1 -a') /  \eta|^2 )$, where $\eta : =  |a' - \tilde{a}|/3$.
Then, $\delta_{a'} ( y^{\ve}_1 ) =   \psi (| ( y^{\ve}_1 -a') /  \eta|^2 ) \cdot \delta_{a'} ( y^{\ve}_1 )$
in ${\bf D}_{- \infty}$.
By Sugita's theorem,  $\theta^{\ve}_{a,a'} (dw)
 =  \psi (| ( a+\Phi_{\ve} ( \tau_{\lambda} (\ve W))^1_{0,1}   -a') /  \eta|^2 )  \cdot \theta^{\ve}_{a,a'} (dw) $,
since $a+\Phi_{\ve} ( \tau_{\lambda} (\ve W))^1_{0,1} $ is the $\infty$-quasi-redefinition 
of $y^{\ve}_1 = y^{\ve} (1, a, w)$.
Hence, 
\begin{align}
\mu_{a,a'}^{\ve} ( B_{R} (X) ) 
&= \theta_{a,a'}^{\ve} (\{ w \in {\cal W} ~|~\ve W \in B_{R} (X) \} )
\nn\\
&=
\int_{  {\cal W} } 
I_{ \{ \ve W \in B_{R} (X) \} }    \cdot 
  \psi (| ( a+\Phi_{\ve} ( \tau_{\lambda} (\ve W))^1_{0,1}   -a') /  \eta|^2 )   \theta^{\ve}_{a,a'} (dw)
=0.
\nn
\end{align}

Next we consider the case $a+ \Phi_{0} ( \tau_{\lambda} (X))^1_{0,1} = a'$.
Note that $\| \ve W^1 - X^1 \|^{4m}_{\al, 4m -B}$ and 
$\| \ve^2 W^2 - X^2 \|^{2m}_{2\al, 2m -B}$ are in $\hat{\cal C}_{4m}$,
even if $X \notin S_2 ({\cal H})$.
Since their $L^2$-norms are bounded in $\ve$, 
so are their $(q,r)$-Sobolev norms for any $(q,r)$.
Note also that $\| D^r y^{\ve}_1\|_q$ is bounded in $\ve$
 for any $(q,r)$.
Set $G(u_1, \ldots, u_n) = \prod_{j=1}^n (u_j - a'_j)^+$.
This is a continuous function from ${\bf R}^n$ to ${\bf R}$
with polynomial growth
and satisfies $\partial_1^2 \cdots \partial_n^2 G (u)=\delta_{a'} (u)$ in 
the sense of Schwartz distributions.

Then, we have 
\begin{align}
\mu_{a,a'}^{\ve} ( B_{R} (X) ) 
&= \theta_{a,a'}^{\ve} (\{ w \in {\cal W} ~|~\ve W \in B_{R} (X) \} )
\nn\\
&\le
{\mathbb E}
\Bigl[
\prod_{i=1}^2
\psi \bigl( \| \ve^i W^i - X^i \|^{4m/i}_{i\al, 4m/i -B} / R^{4m} \bigr)
\cdot
( \partial_1^2 \cdots \partial_n^2 G) (y^{\ve}_1)
\Bigr].
\label{ups1.eq}
\end{align}

Now we use integration by parts as in (\ref{ipb1.eq})--(\ref{ipb2.eq}).
Then, the right hand side of (\ref{ups1.eq}) is equal to 
a finite sum of the following form;
\begin{align}
\sum_{j,k}
{\mathbb E}
\Bigl[
F_{j,k}^{\ve} \cdot
\psi^{(j)} \bigl( \| \ve W^1 - X^1 \|^{4m}_{\al, 4m -B} / R^{4m} \bigr)
\psi^{(k)} \bigl( \| \ve^2 W^2 - X^2 \|^{2m}_{2\al, 2m -B} /  R^{4m} \bigr)
 G (y^{\ve}_1)
\Bigr].
\label{ups2.eq}
\end{align}
Here, $F_{j,k}^{\ve} (w) =F_{j,k}(\ve, w)$ is a polynomial 
in components of 
(i)~$y^{\ve}_1$ and its derivatives,
(ii)~$\| \ve^i W^i - X^i \|^{4m/i}_{i\al, 4m/i -B}$ and its derivatives,
(iii)~ $\tau(\ve)$, which is a Malliavin covariance matrix of $y^{\ve}_1$,
and (iv)~ $\gamma(\ve) = \tau(\ve)^{-1}$.
Note that derivatives of $\gamma(\ve)$ do not appear.

By the uniform non-degeneracy of the Malliavin covariance matrix of $(y^{\ve}_1 -a )/\ve$,
there exists $l >0$ such that $\gamma(\ve)$ is $O(\ve^{-l})$ in any $L^q$.
In turn, this implies that 
there exists $l >0$ such that $F_{j,k}^{\ve}$ is $O(\ve^{-l})$ in any $L^q$.
(Here and in what follows, $l$ may change from line to line.)

Take $q,q' \in (1, \infty)$ so that $1/q + 1/q' =1$.
By H\"older's inequality, 
the right hand side of (\ref{ups2.eq}) is dominated by
\begin{align}
C\ve^{-l}
\mu_{a,a'}^{\ve}
\Bigl( \| \ve^i W^i - X^i \|_{i \al, 4m/i -B}^{1/i}   \le 2^{\frac{1}{m}  } R
~(i=1,2)
 \Bigr)^{\frac{1}{q'} }
=
C\ve^{-l} \mu_{a,a'}^{\ve} (\ve W \in B_{2^{ \frac{1}{m}  }   R} (X))^{\frac{1}{q'} }.
\nn
\end{align}
Here, $C=C_{q,q'} >0$ is a constant independent of $\ve$.
By the large deviation principle of Schilder-type for the law of $\ve W$ on 
$G\Omega^B_{\al, 4m} ({\bf R}^d)$,
\[
\limsup_{\ve \searrow 0} \ve^2 \log \mu_{a,a'}^{\ve} ( B_{R} (X) )  
\le
- \frac{1}{q'}  \inf\{ \|h\|^2_{{\cal H} } /2 ~|~ h \in {\cal H}, S_2(h) \in B_{2^{1/m} R} (X)\}.
\]
Letting $q' \searrow 1$, we have
\begin{eqnarray*}
\limsup_{\ve \searrow 0} \ve^2 \log \mu_{a,a'}^{\ve} ( B_{R} (X) )  
&\le&
- \inf\{ \|h\|^2_{{\cal H} } /2 ~|~ h \in {\cal H}, S_2(h) \in B_{2^{1/m} R} (X)\}
\\
&=&
- \inf\{ I_{Sch}(Y) ~|~ Y\in B_{2^{1/m} R} (X) \}.
\end{eqnarray*}
Since the rate function $I_{Sch}: G\Omega^B_{\al, 4m} ({\bf R}^d) \to [0, \infty]$
is lower semicontinuous, 
the limit of the right hand side as $R \searrow 0$ is dominated by $- I_{Sch}(X)$.
This proves (\ref{limsmball.ineq}).
(Here, $I_{Sch}(S_2(h)) = \|h\|^2_{{\cal H} } /2$ and 
$I_{Sch}(X) =0$ if $X$ is not lying above any $h \in {\cal H}$.)

%
%%%%%%%%%%%%%%%%%%%%%
{\bf [Step 2]}~
The second step is to prove the upper bound in Theorem \ref{tm.ldp.sgt} {\rm (i)}
when $A$ is a compact set.
Let $N \in {\bf N}_+$.
For any $Y \in A$,  take $R =R_{N,Y} >0$ small enough so that 
\[
\limsup_{\ve \searrow 0}  \ve^2 \log \mu_{a,a'}^{\ve} ( B_R (Y) ) \le 
\begin{cases}
-N & ( \mbox{if  $I(Y)= \infty$}), \\
- I(Y) +N^{-1}  &   ( \mbox{if  $I(Y) < \infty$}).
\end{cases}
\]
%
% 
%$\limsup_{\ve \searrow 0}  \ve^2 \log \mu_{a,a'}^{\ve} ( B_r (Y) ) \le - I(Y) +N^{-1}$ if $I(Y) < \infty$.
%
%
The union of such open balls over $Y \in A$ covers the compact set $A$.
Hence, there are finitely many $Y_1, \ldots, Y_k$ such that 
$A \subset  \cup_{i=1}^k B_{R_i} (Y_i)$, where $R_i = R(N, Y_i)$.
By using the remark in the beginning of this section,  we see that
\begin{align}
\limsup_{\ve \searrow 0}  \ve^2 \log \mu_{a,a'}^{\ve} ( A ) 
&\le
(-N ) \vee \max\{ -I(Y_i) +N^{-1}  ~|~ 1 \le i  \le k, \,\,  I(Y_i) < \infty  \}
\nn\\
&\le 
(-N) \vee \bigl(  - \inf_{h \in A \cap {\cal K}^{a,a'}}  \| h\|^2_{{\cal H}}/2   +N^{-1} \bigr).
\nn
\end{align}
Letting $N \to \infty$, we obtain 
\[
\limsup_{\ve \searrow 0}  \ve^2 \log \mu_{a,a'}^{\ve} ( A ) 
\le
- \inf
\{
\| h\|^2_{{\cal H}}/2
~|~ 
h \in {\cal K}^{a,a'}, \, S_2 (h) \in A
\}.
\]
Thus, we have obtained the upper estimate for the compact case.

%%%%%%%%%%%%%%%%%%%%%
{\bf [Step 3]}~
In the final step, we verify the case when $A$ is a closed set. 
Take $\al' \in (\al, 1/2)$ so that $(\al', m)$ still satisfies (\ref{eq.amam}).
As in the proof of Proposition \ref{pr.lowball},
let $B'_R (H)= T_h (B_R)$ be the "ball" of radius $R>0$
centered at $H$ in  $(\al', 4m)$-Besov  geometric rough path space.
(We will use the same symbol by abusing notations.)
Note that this is precompact in $G\Omega^B_{\al, 4m}  ({\bf R}^d)$ 
by Proposition \ref{inj.cpt.pr}.

Choose any $\bar{h} \in {\cal K}^{a,a'} \cap {\cal W}^*$.
For  sufficiently large $R>0$, 
\begin{align}
\mu^{\ve}_{a,a'} ( B'_R ( \bar{H} )^c ) 
&=
\int_{{\cal W}}    
% I_{  \{ \|(\ve W)^i -  \bar{H}^i \|_{i \al' , 4m/i -B}  < R^{i}\, \, (i=1,2)   \}^c }   
I_{  T_{\bar{h}}(B_R)^c } (\ve W )
 \theta^{\ve}_{a,a'}  (dw)
\nn\\
&=
\ve^{-n} 
\int_{{\cal W}}  
\exp \bigl( - \frac{ \la \bar{h} , w \ra}{\ve}  -    \frac{ \|  \bar{h} \|^2_{{\cal H}}  }{2\ve^2} \bigr)
   I_{ \{ \|W^i  \|_{i \al' , 4m/i -B}^{1/i}  <R /\ve
\,\, (i=1,2)\}^c   }   \nu^{\ve} (dw)
\nn\\
&\le
\ve^{-n}   e^{   \| \bar{h}  \|^2_{{\cal H}} /2\ve^2 }
{\mathbb E} \Bigl[
\exp \bigl(  - \frac{ \la 2\bar{h} , w \ra}{\ve}  -    \frac{ \| 2\bar{h} \|^2_{{\cal H}}  }{2\ve^2} \bigr)
\delta_0 \bigl(   \frac{ y^{\ve,  \bar{h} }_1 -a' }{\ve}  \bigr)
\Bigr]^{1/2}
\nn\\
& \qquad \times
\Bigl\{
\sum_{i=1,2}   \nu^{\ve}  ( \{ \|W^i  \|_{i \al' , 4m/i -B}^{1/i}   \ge R /\ve )
\Bigr\}^{1/2}.
\nn
\end{align}

As we proved in (the proof of)
Proposition \ref{pr.lowball} and Lemma \ref{lm.capdcy},  
the last factor above is dominated by 
 $2 e^{ - c_1 R^2 /\ve^2}$ for some constant $c_1 >0$.
Since $\| \delta_{0} ( (y^{\ve, \bar{h} }_1  -a' )/\ve ) \|_{q', -r}$
 is bounded in $\ve$ for some $q', r$, 
it suffices to estimate the $(q,r)$-Sobolev norm 
of the Cameron-Martin density function
$\exp ( -\la 2\bar{h} , w \ra/\ve   -    \| 2\bar{h}  \|^2_{{\cal H}} /2\ve^2 )$.
Using H\"older's inequality,
we can easily see that its $L^q$-norm is dominated by $e^{ 2(q-1) \| \bar{h} \|^2_{{\cal H}}  /\ve^2 }$.
We can also estimate $L^q$-norms of its derivatives and, by using Meyer's equivalence, 
we have 
\[
\Bigl\|  \exp \bigl( \frac{ \la 2\bar{h} , w \ra}{\ve}  -    \frac{ \| 2  \bar{h} \|^2_{{\cal H}}  }{2\ve^2} \bigr) 
\Bigr\|_{q, r} 
 \le
c_2 \ve^{-c_3} (1 +  \| h\|_{{\cal H} })^{c_4}  e^{ 2(q-1) \| \bar{h} \|^2_{{\cal H}}  /\ve^2 }
\]
for some positive constants $c_j =c_{j, q,r}~(j=2,3,4)$, which are independent of $\ve$.
Therefore, we have
\[
\limsup_{ \ve \searrow 0}
\ve^2 \log
\mu^{\ve}_{a,a'} ( B'_R (\bar{H})^c )   \le  - c_1 R^2 + (q -\frac12) \| \bar{h} \|^2_{{\cal H}}.
\]
Note that the right hand side can be made arbitrarily small if we take $R$ large enough.

Now let $A$ be a closed set
such that $A \cap S_2(  {\cal K}^{a,a'} ) \neq \emptyset$.
Take $R>0$ so that 
\[
- c_1 R^2 + (q -\frac12) \| \bar{h}\|^2_{{\cal H}} < - \inf_{h \in A \cap {\cal K}^{a,a'} }  \| h\|^2_{{\cal H}}/2.
\]
Decompose $A = \overline{ B'_R (\bar{H})} \cap A   + ( \overline{ B'_R (\bar{H})} )^c \cap A$
into a union of two disjoint sets.
The second set is included in $B'_R (\bar{H})^c$.
The first set in the disjoint union is compact and, due to the previous step,  we have
\begin{align}
\limsup_{ \ve \searrow 0}
\ve^2 \log
\mu^{\ve}_{a,a'} (  \overline{ B'_R (\bar{H})} \cap A  )   
&\le 
- \inf
\{
\| h\|^2_{{\cal H}}/2
~|~
h \in {\cal K}^{a,a'}, \,  S_2(h) \in A \cap  \overline{ B'_R (\bar{H})} 
\}
\nn\\
&\le
- \inf
\{ \| h\|^2_{{\cal H}}/2
~|~
h \in {\cal K}^{a,a'}, \,  S_2(h) \in A
\}.
\nn
\end{align}
By the remark in the beginning of this section, 
we obtain the upper estimate in this case.
The case $A \cap S_2( {\cal K}^{a,a'}) = \emptyset$
can be done in a similar way and is actually easier.
This completes the proof of Theorem \ref{tm.ldp.sgt}, {\rm (i)}.

%%%%%%%%%%%%%%%%%%%%%%%%%%%%%%%%%%%%%%%%%%%\newpage

\end{document}